\documentclass[11pt]{article}
\usepackage[utf8]{inputenc}
\usepackage{amsthm,amsmath,amssymb,amsfonts}
\usepackage{color}
\usepackage{a4wide}
\usepackage{cancel}
\newcommand{\pder}[2]{\frac{\partial #1}{\partial #2}}

\newcommand{\dx}{\,{\rm d}x}

\newcommand{\Div}{{\rm div}\,}
\newtheorem{thm}{Theorem}
\newtheorem{lem}[thm]{Lemma}

\newtheorem{df}{Definition}

\newcommand{\vr}{\varrho}

\newcommand{\vrd}{\vr_\delta}

\newcommand{\vtd}{\vt_\delta}
\newcommand{\vud}{\vu_\delta}
\newcommand{\vt}{\vartheta}

\newcommand{\vu}{\vc{u}}

\newcommand{\vc}[1]{{\bf #1}}
\newcommand{\vcg}[1]{{\pmb #1}}
\newcommand{\F}[1]{$\mathbb{#1}$}

\newcommand{\tn}[1]{\mbox {\F #1}}

\newcommand{\dS}{\,{\rm d}S}

\newcommand{\R}{\mathbb{R}}

\newcommand{\N}{\mathbb{N}}

\renewcommand{\vrd}{\vr_\delta}
\renewcommand{\vtd}{\vt_\delta}
\renewcommand{\vud}{\vu_\delta}

\begin{document}

\title{Steady compressible Navier--Stokes--Fourier system with general temperature dependent viscosities and hard sphere pressure law}
\author{Zhengguang Guo$^1$ \& Milan Pokorn\' y$^2$}
\maketitle

\centerline{ $^{1}$ School of Mathematics and Statistics}
\centerline{Huaiyin Normal University, Huai'an 223300, Jiangsu, PR China}
\centerline{E-mail: {\tt gzgmath@hytc.edu.cn}}

\centerline{$^{2}$ Charles University, Faculty of Mathematics and Physics}
\centerline{Mathematical Inst. of Charles University, Sokolovsk\' a 83, 186 75 Prague 8, Czech Republic}
\centerline{E-mail: {\tt pokorny@karlin.mff.cuni.cz}}
\vskip0.25cm

\begin{abstract}
We study the existence theory for steady compressible Navier-Stokes-Fourier system in a three dimensional bounded domain for the case of viscosities depending on the temperature in the form $(1+\vartheta)^\alpha$ for $0\leq \alpha\leq 1$ and the hard sphere pressure $p(\vr,\vt) = \vartheta\vr h(\vr)$, $h(\vr)$ is increasing and singular at $\vr=a>0$.  This paper considers both the heat-flux and Dirichlet boundary conditions for the temperature and Dirichlet boundary conditions for the velocity. The key point to controlling the pressure is based on some estimates of the Bogovskii operator.
\end{abstract}


\smallskip


\section{Introduction}
\label{1}

We consider the following steady compressible Navier-Stokes-Fourier system in a bounded domain $\Omega \subset \R^3$ with a sufficiently regular boundary ($\Omega \in C^{1,1}$ at least)
\begin{equation} \label{NSF_1}
\begin{aligned}
\Div(\vr\vu) &=0 \\
\Div(\vr\vu\otimes\vu) -\Div \tn{S} + \nabla p &= \vr \vc{f} \\ 
\Div(\vr\vu E) +\Div\vc{q} -\Div (\tn{S} \vu) + \Div(p\vu) &= \vr \vc{f}\cdot \vu + \vr G.
\end{aligned}
\end{equation}
Here, $\varrho\geq 0$ is the density of the fluid, $\vu$ is the velocity field, $\tn{S}$ is the viscous part of the stress tensor, the scalar function $p$ is the pressure of the fluid which depends on $\varrho$ and thermodynamic temperature $\vartheta >0$, the vector $\vc{q}$ is the heat flux and the scalar $E$ is the specific total energy, $\vc{f}$ and $G$ represent the external force and external heat source, respectively. System (\ref{NSF_1}) describes the steady flow of the heat-conducting compressible fluid, the problem is closed by the set of boundary conditions. In this paper, we  assume that $\vu =\vu_0$ on the boundary and $\vc{q}\cdot\vc{n} = L(\vt-\vt_0)$ or $\vt=\vt_D$. We further distinguish between $\vc{u}_0=\vc{0}$ or $\vc{u}_0$ nonzero. However, in the latter, we always assume that $\vc{u}_0 \cdot \vc{n} =0$, i.e., there is no inflow/outflow of the fluid through the boundary. These conditions will be described in more detail later.

We also need to prescribe the total mass
\begin{equation} \label{mass}
\int_\Omega \vr \dx = M >0.
\end{equation}
Note that the value $M$ cannot be arbitrary; see \eqref{M_limit}.

We consider the stress tensor in the form
\begin{equation} \label{stress}
\tn{S} = \tn{S}(\vt,\nabla\vu) = \mu(\vt) \Big(\nabla \vu + (\nabla \vu)^T-\frac 23 \Div \vu \tn{I}\Big) + \xi (\vt) \Div \vu \tn{I}. 
\end{equation}
The heat flux obeys {\it Fourier's law}
\begin{equation} \label{heat_flow}
\vc{q}=\vc{q}(\vt,\nabla \vt) = -\kappa(\vt)\nabla \vt.
\end{equation}
Here, $(\cdot)^T$ and $\tn{I}$ denote the transpose of a tensor and the identity matrix, respectively.
We assume that the viscosity coefficients behave in the following way. The shear viscosity
$\mu(\vt) >0$ is a Lipschitz continuous function, the bulk viscosity $\xi(\vt) \geq 0$ is a continuous function, and there exist constants $C_1$ and $C_2$ such that
\begin{equation} \label{visc}
\begin{aligned}
C_1 (1+\vt)^\alpha &\leq \mu(\vt) \leq C_2(1+\vt)^\alpha \\
0 &\leq \xi(\vt) \leq C_2(1+\vt)^\alpha, 
\end{aligned}
\end{equation}
where the exponent $\alpha \in [0,1]$. The heat conductivity $\kappa(\vt)>0$ is a continuous function, and there are positive constants $C_3$ and $C_4$ such that 
\begin{equation} \label{heat_cond}
C_3 (1+\vt)^m \leq \kappa (\vt) \leq C_4 (1+\vt)^m
\end{equation}
for a parameter $m>0$. We shall discuss below the choice of the parameters $\alpha$ and $m$ so that we obtain either weak or variational entropy solutions introduced below.

We consider the pressure in the form
\begin{equation} \label{pressure}
p(\vr,\vt)= \vt \vr h(\vr),
\end{equation}
where the function $h(\vr)$ is  non-decreasing in $[0,a)$ with $h(0)=h_0>0$, continuous in $[0,a)$ and continuously differentiable in $(0,a)$ (hence $ h'(\vr)\geq 0$ here) such that
$$
\lim_{z\to a_-} h(z) = \infty.
$$
In fact, all our results remain the same provided we consider 
\begin{equation} \label{pressure_2}
p(\vr,\vt) = \vr h(\vr) + \vr\vt.
\end{equation}
Only the forms of the internal energy and the entropy are slightly different. Note that in order to keep $\varrho < a$, we require
\begin{equation} \label{M_limit}
M < a |\Omega|.
\end{equation}
    
The specific total energy reads
\begin{equation} \label{energy_total}
E(\vr,\vu,\vt) = \frac 12 |\vu|^2 +	e(\vr,\vt),
\end{equation}
where the symbol $e$ stands for the specific internal energy. Using the fact that the pressure and the internal energy are related by the Gibbs relation
\begin{equation} \label{Gibbs}
\frac 1 \vt \Big(De(\vr,\vt) + p(\vr,\vt)D\Big(\frac 1\vr\Big)\Big) = Ds(\vr,\vt),
\end{equation}
where $s(\varrho ,\vartheta )$ is the specific entropy, the pressure $%
p(\varrho ,\vartheta )$ and the internal energy $e(\varrho ,\vartheta )$ are
then related by the Maxwell relation%
\begin{equation}
\frac{\partial e(\varrho ,\vartheta )}{\partial \varrho }=\frac{1}{\varrho
^{2}}\left( p(\varrho ,\vartheta )-\vartheta \frac{\partial p(\varrho
,\vartheta )}{\partial \vartheta }\right) .  \label{maxwell relation}
\end{equation}%
The relative entropy satisfies at least formally the entropy equality
\begin{equation}\label{entropy_equal}
\Div\Big(\frac{\vc{q}}{\vt}\Big) + \Div(\vr s \vu) =\frac{\tn{S}:\nabla \vu }{\vt} - \frac{\vc{q}\cdot\nabla \vt }{\vt^2} + \frac{\vr G}{\vt} 
\end{equation}
Consequently, by (\ref{pressure}) we get 
\begin{equation} \label{energy_internal_general}	
e(\vt,\vr) = g(\vt).
\end{equation}
Often, the relation
\begin{equation} \label{energy_internal_1}	
g(\vt) = c_v \vt
\end{equation}
is used and this form can be used in the case of the heat flux boundary condition. However, in the case of the Dirichlet boundary condition for the temperature, this form is not suitable. The heat capacity at constant volume typically drops to zero as the temperature approaches absolute zero. For example, an ideal Bose gas has $c_v$ proportional to $\vartheta^{3/2}$, and the specific heat of an ideal Fermi gas is proportional to $\vartheta$. Moreover, heat capacities at constant volume and constant pressure approach a common value as $\vartheta\rightarrow 0_+$, see \cite{LL}.   Thus, in some cases, we assume that instead of (\ref{energy_internal_1})
\begin{equation} \label{energy_internal_2}
g(\vt) = \left\{ \begin{array}{rl} 
c_v \vt \quad & \text{ for } \vt >  A^{-\frac 1b} \\
c_v A \vt^{1+b} \quad & \text { for } 0\leq \vt \leq  A^{-\frac 1b}
\end{array}
\right.
\end{equation}
with some $A$, $b>0$. If it is preferable, we may mollify this function in some small neighborhood of the value $\vt=A^{-\frac 1b}$. By suitable choice of $b$ and $A$ we may achieve that this form of the internal energy differs from the previous one only in the set $\vt \in (0,1)$. The corresponding form of the specific entropy is then determined by
\begin{equation} \label{entropy_1}
s(\vr,\vt) = c_v \ln \vt - \int_{\vr_M}^\vr \frac{h(z)}{z} \,{\rm d} z
\end{equation} 
for the case \eqref{energy_internal_1} and 
\begin{equation} \label{entropy_2}
s(\vr,\vt) = \int_{\vt_M}^\vt \frac{g'(s)}{s} \,{\rm d}s - \int_{\vr_M}^\vr \frac{h(z)}{z} \,{\rm d} z
\end{equation}
for the case \eqref{energy_internal_2}. Note that in the latter, for large $\vt$ the first integral grows logarithmically, while for $\vt$ positive, close to zero, it behaves like $\vt^b$ and thus remains bounded for $\vt \to 0_+$.

The first result for the steady compressible Navier--Stokes--Fourier equations for large data goes back to P.L. Lions. The author in \cite{lions} assumed a priori that $\varrho$ is bounded in $L^p(\Omega)$ for $p$ sufficiently large. The heat conducting fluid with only $\varrho\in L^1(\Omega)$  a priori was studied for the first time in \cite{MP1} for $p(\vr,\vt)=\vr^\gamma +\vt \vr$ with $\gamma >3$ with Navier slip boundary conditions for velocity. Then, the authors in \cite{MP2} considered case $\gamma > \frac{7}{3}$ assuming either the slip or no-slip boundary condition for the velocity and the Newton boundary condition for the temperature. All three papers deal with the situation that viscosities are constant and the value of $\gamma$ is far beyond the physically reasonable cases. The case $\gamma>\frac{3}{2}$ which includes a physically realistic case $\gamma=\frac{5}{3}$ was studied in \cite{NP1}. The existence of 
variational entropy solutions was shown in \cite{m3as2014} for any $\gamma >1$ with slip boundary condition for the velocity; moreover, if $\gamma >\frac{5}{4}$,  the solution was shown to be weak. The Dirichlet boundary condition for temperature was considered in \cite{P}. These results were extended to the  general dependence of the viscosities on the
temperature in the form $\mu(\vt),\xi(\vt)\sim(1+\vt)^\alpha$ for $0\leq \alpha \leq 1$ in \cite{Phys}.

In this paper, we focus on the equation of state in general form (\ref{pressure}), the so-called hard sphere pressure. In a certain sense, this formally corresponds to the value $\gamma =\infty$. In particular, this pressure vanishes for $\vt\to 0$, which is particularly relevant to gasses. It is not necessary to add a  ‘‘cold pressure" component independent of $\vt$. On the other hand, we consider very general transport coefficients (\ref{visc}) and (\ref{heat_cond}), where $m$ and $\alpha$ will be related in the existence theory. The hard sphere pressure law also attracted attention in \cite{CFJP}, but in the case of a two dimensional bounded domain. The Trudinger-Moser inequality played an important role in the analysis therein. We present the extension to three dimensional situation including several types of boundary conditions.

\subsection{Formulation of different boundary conditions}

We consider the problem formulated above for different boundary conditions.

 \smallskip
 
{\it Problem 1}:
We consider \eqref{NSF_1}--\eqref{entropy_2} with internal energy either \eqref{energy_internal_1} or \eqref{energy_internal_2} and the corresponding entropies, for the boundary conditions on $\partial \Omega$
\begin{equation} \label{cond_1}
\vu = \vc{0}, \qquad \vc{q}\cdot\vc{n} = L(\vt-\vt_0)
\end{equation}
for a strictly positive function $\vt_0$.

\bigskip

{\it Problem 2}:
We consider \eqref{NSF_1}--\eqref{entropy_2} with internal energy \eqref{energy_internal_2} and the corresponding entropy \eqref{entropy_2}, for the boundary conditions on $\partial \Omega$
\begin{equation} \label{cond_2}
\vu = \vc{0}, \qquad \vt=\vt_D
\end{equation}
for a strictly positive $\vt_D$ with
$\vartheta _{D}\in W^{2-\frac{1}{q},q}(\partial \Omega )$ for some $q>3$.
\bigskip

{\it Problem 3}:
We consider \eqref{NSF_1}--\eqref{entropy_2} with internal energy either \eqref{energy_internal_1} or \eqref{energy_internal_2} and the corresponding entropies, for the boundary conditions on $\partial \Omega$
\begin{equation} \label{cond_3}
\vu = \vu_{0}, \qquad \vc{q}\cdot\vc{n} = L(\vt-\vt_0)
\end{equation}
for a strictly positive function $\vt_0$ and $\vu_0$ non-zero. This boundary data  $\vu_{0}$ is assumed to be extended to the whole $\Omega$ so that
\begin{equation} \label{cond_3a}
    \vu_{0}\in W^{1,\infty}(\Omega), \quad \text{div}\, \vu_{0}=0, \quad \vu_0 \cdot \vc{n}|_{\partial \Omega} =0.
\end{equation}

\bigskip

{\it Problem 4}:
We consider \eqref{NSF_1}--\eqref{entropy_2} with internal energy \eqref{energy_internal_2} and the corresponding entropy \eqref{entropy_2}, for the boundary conditions on $\partial \Omega$
\begin{equation} \label{cond_4}
\vu = \vu_{0}, \qquad \vt=\vt_D
\end{equation}
for a strictly positive 
$\vartheta _{D}\in W^{2-\frac{1}{q},q}(\partial \Omega )$ for some $q>3$, and non-zero $\vu_0\in  W^{1,\infty}(\Omega)$ with $ \vu_0 \cdot \vc{n}|_{\partial \Omega} =0$.

\section{Weak formulation, main results}

Our definition of weak and variational solutions is based on the following identities and inequalities.

\medskip

Weak form of the continuity equation:
\begin{equation} \label{weak_cont}
\int_\Omega \vr \vu \cdot \nabla \psi \dx = 0
\end{equation}
for all $\psi \in C^1(\overline{\Omega})$.

\medskip

Weak form of the momentum equation:
\begin{equation} \label{weak_mom}
\int_\Omega \Big[ \vr (\vu \otimes \vu):\nabla \vcg{\varphi} + p(\vr,\vt) \Div \vcg{\varphi} -\tn{S}(\vt,\nabla \vu) : \nabla \vcg{\varphi} + \vr\vc{f}\cdot \vcg{\varphi}\Big]\dx = 0
\end{equation}
for all $\vcg{\varphi}\in C^1_0(\Omega;\R^3)$.

\medskip

We have to distinguish different situations for the total energy balance. For the zero velocity on the boundary and heat flux boundary condition, we use the following definition. \newline 
\begin{equation} \label{weak_energy_1}
\begin{aligned}
&\int_\Omega \Big[\vr \Big(\frac 12 |\vu|^2 + e(\vr,\vt)\Big)\vu \cdot \nabla \psi + p(\vr,\vt)\vu \cdot \nabla \psi -\tn{S}(\vt,\nabla \vu)\vu \cdot \nabla \psi +\vc{q}(\vt,\nabla \vt) \cdot \nabla \psi \Big]\dx \\
& + \int_\Omega \Big(\vr \vc{f}\cdot \vu + \vr G\Big)\psi \dx = \int_{\partial \Omega} L(\vt-\vt_0)\psi\dS 
\end{aligned}
\end{equation}
for all $\psi \in C^1(\overline{\Omega})$.

\smallskip

For the zero velocity on the boundary, but Dirichlet boundary condition for the temperature, we consider
\begin{equation} \label{weak_energy_2}
\begin{aligned}
&\int_\Omega \Big[\vr \Big(\frac 12 |\vu|^2 + e(\vr,\vt)\Big)\vu \cdot \nabla \psi + p(\vr,\vt)\vu \cdot \nabla \psi -\tn{S}(\vt,\nabla \vu)\vu \cdot \nabla \psi +\vc{q}(\vt,\nabla \vt) \cdot \nabla \psi \Big]\dx \\
& + \int_\Omega \Big(\vr \vc{f}\cdot  \vu + \vr G\Big)\psi \dx = 0 
\end{aligned}
\end{equation}
for all $\psi \in C^1_0(\Omega)$.

For the non-zero velocity (but zero normal trace) on the boundary and the heat flux boundary condition, we need to combine the momentum equation (however, formally with test function $\psi \vu_0$) and the ‘‘total energy balance" with test function $\psi$, $\psi \in C^1(\overline{\Omega})$. It reads
\begin{equation} \label{weak_energy_3}
\begin{aligned}
&\int_\Omega \Big[\vr \Big(\frac 12 |\vu|^2 + e(\vr,\vt)\Big)\vu \cdot \nabla \psi + p(\vr,\vt)\vu \cdot \nabla \psi -\tn{S}(\vt,\nabla \vu)\vu \cdot \nabla \psi +\vc{q}(\vt,\nabla \vt) \cdot \nabla \psi  \Big]\dx \\
&  + \int_\Omega  \big(\vr \vc{f}\cdot  \vu + \vr G\big)\psi \dx -\int_{\partial \Omega} L(\vt-\vt_0)\psi\dS\\
&=  \int_\Omega \Big[ \vr (\vu \otimes \vu):\nabla (\psi \vu_0) + p(\vr,\vt) \Div (\psi \vu_0) -\tn{S}(\vt,\nabla \vu) : \nabla (\psi\vu_0) + \vr\vc{f}\cdot \vu_0 \psi \Big]\dx 
\end{aligned}
\end{equation}
for all $\psi \in C^1(\overline{\Omega})$.

We also distinguish several variants for the entropy inequality. For the zero velocity on the boundary and heat flux boundary condition, we consider

\begin{equation} \label{weak_entropy_1}
\begin{aligned}
&\int_\Omega \Big(\frac{\tn{S}(\vt,\nabla \vu):\nabla \vu }{\vt} - \frac{\vc{q}(\vt,\nabla \vt)\cdot\nabla \vt }{\vt^2} + \frac{\vr G}{\vt}\Big)\psi \dx \\
&\leq \int_{\partial \Omega} \frac{L(\vt-\vt_0)}{\vt}\psi \dS -\int_\Omega \Big(\frac{\vc{q}(\vt,\nabla \vt)\cdot\nabla \psi}{\vt}+ \vr s(\vr,\vt)\vu \cdot \nabla \psi\Big)\dx
\end{aligned} 
\end{equation}
for all non-negative $\psi \in C^1(\overline{\Omega})$. 

For the  other cases, we only assume that

\begin{equation} \label{weak_entropy_2}
\begin{aligned}
&\int_\Omega \Big(\frac{\tn{S}(\vt,\nabla \vu):\nabla \vu }{\vt} - \frac{\vc{q}(\vt,\nabla \vt)\cdot\nabla \vt }{\vt^2} + \frac{\vr G}{\vt}\Big)\psi \dx \\
&\leq -\int_\Omega \Big(\frac{\vc{q}(\vt,\nabla \vt)\cdot\nabla \psi}{\vt}+ \vr s(\vr,\vt)\vu \cdot \nabla \psi\Big)\dx
\end{aligned} 
\end{equation}
for all $\psi \in C^1_0(\Omega)$, non-negative.

For zero velocity on the boundary and Dirichlet boundary condition for the temperature, we consider the ballistic energy inequality in the form

\begin{equation} \label{ballistic_energy_1}
\begin{aligned}
&\int_\Omega \Big(\frac{\tn{S}(\vt,\nabla \vu):\nabla \vu }{\vt} - \frac{\vc{q}(\vt,\nabla \vt):\nabla \vt }{\vt^2} + \frac{\vr G}{\vt}\Big)\widetilde{\vt} \dx \\
&\leq -\int_\Omega \Big(\frac{\vc{q}(\vt,\nabla \vt)\cdot\nabla \widetilde{\vt}}{\vt}+ \vr s(\vr,\vt)\vu \cdot \nabla \widetilde{\vt}\Big)\dx
+\int_\Omega \big(\vr G + \vr \vc{f}\cdot \vu\big) \dx 
\end{aligned} 
\end{equation}
for any $\widetilde{\vt}$ being an extension of the boundary data $\vt_D$ to $\Omega$.

Finally, for the case of general $\vu = \vu_0$ and $\vt=\vt_D$ on $\partial \Omega$ we require the following modified ballistic energy inequality
\begin{equation} \label{ballistic_energy_2}
\begin{aligned}
&\int_\Omega \Big(\frac{\tn{S}(\vt,\nabla \vu):\nabla \vu }{\vt} - \frac{\vc{q}(\vt,\nabla \vt):\nabla \vt }{\vt^2} + \frac{\vr G}{\vt}\Big)\widetilde{\vt} \dx \\
&\leq -\int_\Omega \Big(\frac{\vc{q}(\vt,\nabla \vt)\cdot\nabla \widetilde{\vt}}{\vt}+ \vr s(\vr,\vt)\vu \cdot \nabla \widetilde{\vt}\Big)\dx
+\int_\Omega  \big(\vr \vc{f}\cdot \vu + \vr G\big) \dx 
 \\ 
&+ \int_\Omega \Big[ \vr (\vu \otimes \vu):\nabla \vu_0 + p(\vr,\vt) \Div  \vu_0 -\tn{S}(\vt,\nabla \vu) : \nabla \vu_0 + \vr\vc{f}\cdot \vu_0  \Big]\dx
\end{aligned} 
\end{equation}
for any $\widetilde{\vt}$ being an extension of the boundary data $\vt_D$ to $\Omega$.

In each definition, we assume that the unknown functions are such that all integrals are finite. Note also that any weak solution to the continuity equation (if $\nabla \vu \in L^r(\Omega)$ for some $r>1$) is immediately also a renormalized solution. We start with Problem 1. 

\begin{df} \label{def_p1}
We say that $(\vr,\vu,\vt)$ is a variational entropy solution to Problem 1 provided we have $\vu \in W^{1,q}_0(\Omega)$ for some $q>1$, $\ln \vt$, $\vt^{\frac m2} \in W^{1,2}(\Omega)$, we have \eqref{weak_cont}, \eqref{weak_mom}, \eqref{weak_entropy_1} as well as \eqref{weak_energy_1} for $\psi \equiv 1$. The variational entropy solution is also a weak solution, provided \eqref{weak_energy_1} holds for arbitrary $\psi \in C^1(\overline{\Omega})$.
\end{df}

Here, we read our estimates from \eqref{weak_entropy_1} and \eqref{weak_energy_1}, both with $\psi \equiv 1$. This guarantees that all terms in other equalities but the pressure are integrable. We only need to estimate the pressure, which will be done in two steps. As $\varrho$ is by the form of the pressure in $L^\infty(\Omega)$, the renormalized form of the continuity equation is a straightforward consequence of the Friedrichs commutator lemma. Finally, we obtain strong convergence of the approximate  sequence for the density using standard tools. Under more restrictive assumptions, we get \eqref{weak_energy_1} for any $\psi \in C^1(\overline{\Omega})$.

Next, we continue with Problem 2
\begin{df} \label{def_p2}
We say that $(\vr,\vu,\vt)$ is a variational entropy solution to Problem 2 provided we have $\vu \in W^{1,q}_0(\Omega)$ for some $q>1$, $\ln \vt$, $\vt^{\frac m2} \in W^{1,2}(\Omega)$, $\vt=\vt_D$ on $\partial \Omega$ in the sense of traces and we have \eqref{weak_cont}, \eqref{weak_mom}, \eqref{weak_entropy_2} as well as \eqref{ballistic_energy_1}. The solution is weak, provided also  \eqref{weak_energy_2} holds.
\end{df}

We read the necessary estimates from \eqref{ballistic_energy_1}. The possibility of estimating the right-hand side of the ballistic energy inequality is closely connected with the form of the entropy, which cannot be singular for temperature close to zero. The rest of the arguments is similar to Problem 1.

For Problem 3, we have

\begin{df} \label{def_p3}
We say that $(\vr,\vu,\vt)$ is a variational entropy solution to Problem 3 provided we have $\vu \in W^{1,q}(\Omega)$ for some $q>1$, $\vu=\vu_0$ at $\partial \Omega$ in the sense of traces, $\ln \vt$, $\vt^{\frac m2} \in W^{1,2}(\Omega)$, and we have \eqref{weak_cont}, \eqref{weak_mom}, \eqref{weak_entropy_1} as well as \eqref{weak_energy_3} for $\psi=1$. The solution is weak if we additionally have \eqref{weak_energy_3} for arbitrary $\psi \in C^1(\overline{\Omega})$.
\end{df}

We now read our estimates from \eqref{weak_energy_3} and \eqref{weak_entropy_1}, both with $\psi \equiv 1$. The form of the entropy for small temperature is not important in this case and the procedure is similar to Problems 1 and 2.

Finally, for Problem 4 we have

\begin{df} \label{def_p4}
We say that $(\vr,\vu,\vt)$ is a weak solution to Problem 4 provided we have $\vu \in W^{1,q}(\Omega)$ for some $q>1$, $\vu=\vu_0$ at $\partial \Omega$ in the sense of traces, $\ln \vt$, $\vt^{\frac m2} \in W^{1,2}(\Omega)$, $\vt=\vt_D$ on $\partial \Omega$ in the sense of traces and we have \eqref{weak_cont}, \eqref{weak_mom}, \eqref{weak_entropy_2}, \eqref{ballistic_energy_2} as well as \eqref{weak_energy_2}.
\end{df}

 We read the estimates from \eqref{ballistic_energy_2} which leads to the requirement of a certain smallness of $\nabla \vc{u}_0$. The entropy needs to remain bounded for temperature close to zero, similarly as in Problem 2. 

\subsection{Existence of solutions}

In the whole paper, we assume that $\Omega \in C^{1,1}$. This makes the construction of approximate solutions possible, even though, for the formal a priori estimates of the solution, it is enough to assume that the domain is Lipschitz. In some cases, it is possible to approximate the Lipschitz domain by smooth ones and obtain the solution also for less regular domains, however, we shall not consider this. Moreover, we also assume that $\vc{f}\in L^\infty(\Omega;\R^3)$ and $G\in L^\infty(\Omega)$, although less regularity is in some cases sufficient. Finally, we assume that $G \geq 0$ a.e. in $\Omega$. The main results are as follows.

\begin{thm} \label{t1}
Under the assumptions stated above, let $\alpha \in [0,1]$ and $m > \max\{\frac{1+\alpha}{3}, \frac{1-\alpha}{2}\}$. Then there exists a variational entropy solution to Problem 1. Furthermore, if $m>1$, then the solution is also weak.
\end{thm}

\begin{thm} \label{t2}
Under the assumptions stated above, let $\alpha \in [0,1]$ and $m > \max\{\frac{1+\alpha}{3}, 1-\alpha\}$. Then there exists a variational entropy solution to Problem 2. Furthermore, if $m>1$, then the solution is also weak.
\end{thm}

\begin{thm} \label{t3}
Under the assumptions stated above, let $\alpha \in (0,1)$ and $m > \max\{\frac{1+\alpha}{3}, \frac{1-\alpha}{2}\}$. Then there exists a variational entropy solution to Problem 3. Furthermore, if $m>1$, then the solution is also weak. Furthermore, if $\alpha =0$ and $m>1$, then there exists a variational entropy solution to Problem 3 that is also weak. Finally, if $\alpha=1$, $\|\nabla \vc{u}_0\|_\infty$ is sufficiently small and $m\geq \frac 23$, then there exists a variational entropy solution to Problem 3 which is for $m>1$ also a weak solution. 
\end{thm}

\begin{thm} \label{t4}
Under the assumptions stated above, let $\alpha =1$, $m\geq 2$, and let  $\|\nabla \vu_0\|_{\infty}$ be sufficiently small. Then there exists a weak entropy solution to Problem 4. 
\end{thm}

Before presenting the proof, we introduce an important generalization of the
Korn inequality, which justifies our choice of the structure of $\mbox {\F S}%
(\vartheta ,\nabla \mathbf{u})$. Its proof can be found in \cite{Feireisl2017} with a zero extension of $\mathbf{u}$ or in \cite{m3as2014} for $\mathbf{u}\in W_{
\mathbf{n}}^{1,p}(\Omega;\mathbb{R}^{3}):=  \{\mathbf{u\in }W^{1,p}(\Omega ;
\mathbb{R}^{3})\colon\mathbf{u\cdot n}=0$ in the sense of traces on $\partial
\Omega \}$.
\begin{lem}
\label{Korn ineq}Let $1<p<\infty $ and let the tensor function $\mbox {\F S}%
(\vartheta ,\nabla\mathbf{u})$ satisfy (\ref{visc}). Then we have either for $\mathbf{u}%
\in W_{0}^{1,p}(\Omega ;\mathbb{R}^{3})$ or for $\mathbf{u}\in W_{%
\mathbf{n}}^{1,p}(\Omega ;\mathbb{R}^{3})$ (in the latter, we additionally assume that $\Omega$ is not axially symmetric) that%
\begin{equation*}
\int_{\Omega }\left\vert \left( \nabla \mathbf{u+}\nabla \mathbf{u}^{T}-%
\frac{2}{3}\mathrm{div}\,\mathbf{u}\mbox {\F I}\right) :\nabla \mathbf{u}%
\right\vert ^{\frac{p}{2}}\dx\geq C\left\Vert \mathbf{u}\right\Vert _{1,p}^{p}%
\text{ }.
\end{equation*}
\end{lem}

\section{Approximation}

We will not discuss in detail the construction of the approximate system since it can be performed very similarly to \cite{CFJP}. The corresponding approximate system reads as follows:

We take a smooth non-decreasing function $T$: $\R \to [-2a,2a]$ such that $T(\vr) = \vr$ for $0\leq \vr\leq a$ and the rest of $\R$ it smoothly maps to $[-2a,2a]$. We fix
\begin{equation} \label{appr_par}
\varepsilon \in \R^+, \quad N \in \N, \quad \eta \in \R^+, \quad R\in \R^+, \quad \delta \in \R^+.
\end{equation}
For fixed $N$ we set
$$
X_N = {\rm span}\{\vc{w}_1,\vc{w}_2,\dots,\vc{w}_N\},
$$
where $\{\vc{w}\}_{i=1}^\infty$ forms a complete orthonormal system in $L^2(\Omega;\R^3)$ and a complete orthogonal system in $W^{1,2}_0(\Omega;\R^3)$ (e.g., eigenfunctions of the vector-valued Laplace operator with homogeneous Dirichlet boundary conditions). We then consider the following approximate system ($\vu -\vu_0 \in X_N$, where $\vu_0 =\vc{0}$ in case we consider the homogeneous Dirichlet boundary conditions for the velocity)
\begin{equation} \label{app_cont}
\begin{aligned}
\Div (\vr\vu) &= \varepsilon \Delta \vr -\varepsilon \Big(\vr-\frac{M}{|\Omega|}\Big) &&\qquad \text{ in } \Omega, \\
\pder{\vr}{\vc{n}} &= 0 &&\qquad \text{ on } \partial \Omega,
\end{aligned}
\end{equation}
\begin{equation}\label{app_moment}
\begin{aligned}
&\int_\Omega \Big(\frac 12 T(\vr) (\vu \otimes \vu:\nabla \vcg{\varphi} - (\vu \cdot \nabla \vu)\cdot \varphi ) + (p_R (\vr,\vt) + \varepsilon (\vr^\beta + \vr^2)) \Div\vcg{\varphi} \Big)\dx \\
& = \int_\Omega \Big(\tn{S}^\eta_\delta (\vt,\nabla \vu) :\nabla \vcg{\varphi} - T(\vr) \vc{f}\cdot \vcg{\varphi}\Big)\dx
\end{aligned}
\end{equation}
for all $\vcg{\varphi} \in X_N$,
\begin{equation} \label{app_energ}
\begin{aligned}
 \Div(\vr e_R(\vt)\vu) + \Div \vc{q}^\eta_\delta(\vt,\nabla \vt) &= \tn{S}^\eta_\delta(\vt,\nabla \vu):\nabla \vu -p_R(\vr,\vt)\Div \vu \\
 & + \varepsilon^2 |\nabla \vr|^2 (\beta \vr^{\beta-2} + 2) +\vr (G+\varepsilon) &&\quad \text{ in } \Omega,  \\
\text{ either } \vc{q}^\eta_\delta \cdot\vc{n} &= L(\vt-\vt_0) && \quad \text{ on } \partial \Omega, \\
\text{ or } \vt & = \vt_D && \quad \text{ on } \partial \Omega.
\end{aligned}
\end{equation}

We use (below, the upper index $\eta$ always refers to the mollification of a function)
\begin{equation} \label{app_stress}
\tn{S}^\eta_\delta (\vt,\nabla \vu) = \frac{\mu^\eta(\vt) + \delta \vt}{\eta\vt +1}\big(\nabla \vu + (\nabla \vu)^T - \frac 23 \Div \vu \tn{I}\big) + \frac{\xi^\eta(\vt) + \delta \vt}{\eta\vt +1} \Div \vu \tn{I}
\end{equation}
\begin{equation} \label{app_heatfl}
\vc{q}^\eta_\delta(\vt,\nabla \vt) = -\Big(\kappa^\eta(\vt) + \delta(\vt^B+ \vt^{-1})\Big) \nabla \vt,  
\end{equation}
where $\beta$, $B \gg 1$. If $\alpha =1$, we do not consider the extra term $\delta \vt$ in the viscosity regularization, as this term is useless, and we write $\tn{S}$ instead of $\tn{S}_\delta$ in what follows. We further set the following
\begin{equation} \label{app_thermo}
\begin{aligned}
p_R(\vr,\vt) = \vr\vt h_R(\vt), \qquad e_R(\vt) = e(\vt)= g(\vt), \\
s_R(\vr,\vt)= \int_{\vt_M}^\vt \frac{g'(z)}{z} \,{\rm d}z - \int_{\vr_M}^\vr \frac{h_R(z)}{z} \,{\rm d}z,
\end{aligned}
\end{equation}
where 
$$
h_R(\vr) = \left\{
\begin{array}{rl} h(\vr), & 0 \leq \vr \leq a-\frac 1R, \\
h\Big(a-\frac 1R\Big) + h'\Big(a-\frac 1R\Big) \Big(\vr-a + \frac 1R\Big), & \text{ otherwise }.
\end{array}
\right.
$$
Note also that the Gibbs relation \eqref{Gibbs} holds for $p_R$, $s_R$ and $e_R$. Note that the only change with respect to the approximation in \cite{CFJP} is the extra regularization of the pressure which also reflects in the internal energy balance. Recall that our constructed sequence of densities is non-negative and of temperatures positive a.e. in $\Omega$ and the same will hold for the limit functions.

We now subsequently send $N\to\infty$, $\eta\to 0_+$, $R\to \infty$ and $\varepsilon \to 0_+$. Since the computations are very similar to the situation in \cite{CFJP}, we skip the details. Note that we need the regularization of the pressure in order to have better information about the density before we pass with $R\to \infty$. After the limit passage, we know that the density is bounded, hence we can remove these terms in the subsequent limit $\varepsilon \to 0_+$. Thus, all computations are, in fact, a combination of the approach in \cite{CFJP} and \cite{NP1}. We end up with the following formulations, in dependence on the problem we want to solve.
In what follows, we denote the sequence of approximate solutions $(\vrd,\vud,\vtd)$.

\bigskip

{\bf Problem 1:}

\begin{equation} \label{weak_cont_delta_1}
\int_\Omega \vrd \vud \cdot \nabla \psi \dx = 0
\end{equation}
for all $\psi \in C^1(\overline{\Omega})$;

\smallskip

\begin{equation} \label{weak_mom_delta_1}
\int_\Omega \Big[ \vrd (\vud \otimes \vud):\nabla \vcg{\varphi} + p(\vrd,\vtd) \Div \vcg{\varphi} -\tn{S}_\delta(\vtd,\nabla \vud) : \nabla \vcg{\varphi} + \vrd\vc{f}\cdot \vcg{\varphi}\Big]\dx = 0
\end{equation}
for all $\vcg{\varphi}\in C^1_0(\Omega;\R^3)$;

\begin{equation} \label{weak_energy_1_delta_1}
\begin{aligned}
&\int_\Omega \Big[\vrd \Big(\frac 12 |\vud|^2 + e(\vtd)\Big)\vud \cdot \nabla \psi + p(\vrd,\vtd)\vud \cdot \nabla \psi -\tn{S}_\delta(\vtd,\nabla \vud)\vud \cdot \nabla \psi +\vc{q}_\delta(\vtd,\nabla \vtd) \cdot \nabla \psi \Big] \\
& + \int_\Omega \Big(\vrd \vc{f}\cdot \vud+ \vrd G\Big)\psi \dx = \int_{\partial \Omega} L(\vtd-\vt_0)\psi\dS 
\end{aligned}
\end{equation}
for all $\psi \in C^1(\overline{\Omega})$ and

\begin{equation} \label{weak_entropy_1_delta_1}
\begin{aligned}
&\int_\Omega \Big(\frac{\tn{S}_\delta(\vtd,\nabla \vud):\nabla \vud }{\vtd} - \frac{\vc{q}_\delta(\vtd,\nabla \vtd)\cdot\nabla \vtd }{\vtd^2} + \frac{\vrd G}{\vtd}\Big)\psi \dx \\
&\leq \int_{\partial \Omega} \frac{L(\vtd-\vt_0)}{\vtd}\psi \dS -\int_\Omega \Big(\frac{\vc{q}_\delta(\vtd,\nabla \vtd)\cdot\nabla \psi}{\vtd}+ \vrd s(\vrd,\vtd)\vud \cdot \nabla \psi\Big)\dx
\end{aligned} 
\end{equation}
for all non-negative $\psi \in C^1(\overline{\Omega})$.

\bigskip

{\bf Problem 2:}

\begin{equation} \label{weak_cont_delta_2}
\int_\Omega \vrd \vud \cdot \nabla \psi \dx = 0
\end{equation}
for all $\psi \in C^1(\overline{\Omega})$;

\smallskip

\begin{equation} \label{weak_mom_delta_2}
\int_\Omega \Big[ \vrd (\vud \otimes \vud):\nabla \vcg{\varphi} + p(\vrd,\vtd) \Div \vcg{\varphi} -\tn{S}_\delta(\vtd,\nabla \vud) : \nabla \vcg{\varphi} + \vrd\vc{f}\cdot \vcg{\varphi}\Big]\dx = 0
\end{equation}
for all $\vcg{\varphi}\in C^1_0(\Omega;\R^3)$;

\begin{equation} \label{weak_energy_2_delta_2}
\begin{aligned}
&\int_\Omega \Big[\vrd \Big(\frac 12 |\vud|^2 + e(\vtd)\Big)\vud \cdot \nabla \psi + p(\vrd,\vtd)\vud \cdot \nabla \psi -\tn{S}_\delta(\vtd,\nabla \vud)\vud \cdot \nabla \psi +\vc{q}_\delta(\vtd,\nabla \vtd) \cdot \nabla \psi \Big] \\
& + \int_\Omega \Big(\vrd \vc{f}\cdot  \vud + \vrd G\Big)\psi \dx = 0 
\end{aligned}
\end{equation}
for all $\psi \in C^1_0(\Omega)$;

\begin{equation} \label{weak_entropy_2_delta_2}
\begin{aligned}
&\int_\Omega \Big(\frac{\tn{S}_\delta(\vtd,\nabla \vud):\nabla \vud }{\vtd} - \frac{\vc{q}_\delta(\vtd,\nabla \vtd)\cdot\nabla \vtd }{\vtd^2} + \frac{\vrd G}{\vtd}\Big)\psi \dx \\
&\leq -\int_\Omega \Big(\frac{\vc{q}_\delta(\vtd,\nabla \vtd)\cdot\nabla \psi}{\vtd}+ \vrd s(\vrd,\vtd)\vud \cdot \nabla \psi\Big)\dx
\end{aligned} 
\end{equation}
for all $\psi \in C^1_0(\Omega)$, non-negative and

\begin{equation} \label{ballistic_energy_1_delta_2}
\begin{aligned}
&\int_\Omega \Big(\frac{\tn{S}_\delta(\vtd,\nabla \vud):\nabla \vud }{\vtd} - \frac{\vc{q}_\delta(\vtd,\nabla \vtd):\nabla \vtd }{\vtd^2} + \frac{\vrd G}{\vtd}\Big)\widetilde{\vt} \dx \\
&\leq -\int_\Omega \Big(\frac{\vc{q}_\delta(\vtd,\nabla \vtd)\cdot\nabla \widetilde{\vt}}{\vtd}+ \vrd s(\vrd,\vtd)\vud \cdot \nabla \widetilde{\vt}\Big)\dx
+\int_\Omega \big(\vrd G + \vrd \vc{f}\cdot \vud\big) \dx 
\end{aligned} 
\end{equation}
for any $\widetilde{\vt}$ being an extension of the boundary data $\vt_D$ to $\Omega$.

\bigskip

{\bf Problem 3:}

\begin{equation} \label{weak_cont_delta_3}
\int_\Omega \vrd \vud \cdot \nabla \psi \dx = 0
\end{equation}
for all $\psi \in C^1(\overline{\Omega})$;

\smallskip

\begin{equation} \label{weak_mom_delta_3}
\int_\Omega \Big[ \vrd (\vud \otimes \vud):\nabla \vcg{\varphi} + p(\vrd,\vtd) \Div \vcg{\varphi} -\tn{S}_\delta(\vtd,\nabla \vud) : \nabla \vcg{\varphi} + \vrd\vc{f}\cdot \vcg{\varphi}\Big]\dx = 0
\end{equation}
for all $\vcg{\varphi}\in C^1_0(\Omega;\R^3)$;

\begin{equation} \label{weak_energy_3_delta_3}
\begin{aligned}
&\int_\Omega \Big[\vrd \Big(\frac 12 |\vud|^2 + e(\vtd)\Big)\vud \cdot \nabla \psi + p(\vrd,\vtd)\vud \cdot \nabla \psi -\tn{S}_\delta(\vtd,\nabla \vud)\vud \cdot \nabla \psi +\vc{q}_\delta(\vtd,\nabla \vtd) \cdot \nabla \psi  \Big] \\
&  + \int_\Omega  \big(\vrd \vc{f}\cdot  \vud + \vrd G\big)\psi \dx -\int_{\partial \Omega} L(\vtd-\vt_0)\psi\dS\\
&=  \int_\Omega \Big[ \vrd (\vud \otimes \vud):\nabla (\psi \vu_0) + p(\vrd,\vtd) \Div (\psi \vu_0) -\tn{S}_\delta(\vtd,\nabla \vud) : \nabla (\psi\vu_0) + \vrd\vc{f}\cdot \vu_0 \psi \Big]\dx 
\end{aligned}
\end{equation}
for all $\psi \in C^1(\overline{\Omega})$ and

\begin{equation} \label{weak_entropy_1_delta_3}
\begin{aligned}
&\int_\Omega \Big(\frac{\tn{S}_\delta(\vtd,\nabla \vud):\nabla \vud }{\vtd} - \frac{\vc{q}_\delta(\vtd,\nabla \vtd)\cdot\nabla \vtd }{\vtd^2} + \frac{\vrd G}{\vtd}\Big)\psi \dx \\
&\leq \int_{\partial \Omega} \frac{L(\vtd-\vt_0)}{\vtd}\psi \dS -\int_\Omega \Big(\frac{\vc{q}_\delta(\vtd,\nabla \vtd)\cdot\nabla \psi}{\vtd}+ \vrd s(\vrd,\vtd)\vud \cdot \nabla \psi\Big)\dx
\end{aligned} 
\end{equation}
for all non-negative $\psi \in C^1(\overline{\Omega})$.

\bigskip

{\bf Problem 4:}

\begin{equation} \label{weak_cont_delta_4}
\int_\Omega \vrd \vud \cdot \nabla \psi \dx = 0
\end{equation}
for all $\psi \in C^1(\overline{\Omega})$;

\smallskip

\begin{equation} \label{weak_mom_delta_4}
\int_\Omega \Big[ \vrd (\vud \otimes \vud):\nabla \vcg{\varphi} + p(\vrd,\vtd) \Div \vcg{\varphi} -\tn{S}_\delta(\vtd,\nabla \vud) : \nabla \vcg{\varphi} + \vrd\vc{f}\cdot \vcg{\varphi}\Big]\dx = 0
\end{equation}
for all $\vcg{\varphi}\in C^1_0(\Omega;\R^3)$;

\begin{equation} \label{weak_energy_2-delta_4}
\begin{aligned}
&\int_\Omega \Big[\vrd \Big(\frac 12 |\vud|^2 + e(\vtd)\Big)\vud \cdot \nabla \psi + p(\vrd,\vtd)\vud \cdot \nabla \psi -\tn{S}_\delta(\vtd,\nabla \vud)\vud \cdot \nabla \psi +\vc{q}_\delta(\vtd,\nabla \vtd) \cdot \nabla \psi \Big] \\
& + \int_\Omega \Big(\vrd \vc{f}\cdot  \vud + \vrd G\Big)\psi \dx = 0 
\end{aligned}
\end{equation}
for all $\psi \in C^1_0(\Omega)$;

\begin{equation} \label{weak_entropy_2_delta_4}
\begin{aligned}
&\int_\Omega \Big(\frac{\tn{S}_\delta(\vtd,\nabla \vud):\nabla \vud }{\vtd} - \frac{\vc{q}_\delta(\vtd,\nabla \vtd)\cdot\nabla \vtd }{\vtd^2} + \frac{\vrd G}{\vtd}\Big)\psi \dx \\
&\leq -\int_\Omega \Big(\frac{\vc{q}_\delta(\vtd,\nabla \vtd)\cdot\nabla \psi}{\vtd}+ \vrd s(\vrd,\vtd)\vud \cdot \nabla \psi\Big)\dx
\end{aligned} 
\end{equation}
for all $\psi \in C^1_0(\Omega)$, non-negative and

\begin{equation} \label{ballistic_energy_2_delta_4}
\begin{aligned}
&\int_\Omega \Big(\frac{\tn{S}_\delta(\vtd,\nabla \vud):\nabla \vud }{\vtd} - \frac{\vc{q}_\delta(\vtd,\nabla \vtd):\nabla \vtd }{\vtd^2} + \frac{\vrd G}{\vtd}\Big)\widetilde{\vt} \dx \\
&\leq -\int_\Omega \Big(\frac{\vc{q}_\delta(\vtd,\nabla \vtd)\cdot\nabla \widetilde{\vt}}{\vtd}+ \vrd s(\vrd,\vtd)\vud \cdot \nabla \widetilde{\vt}\Big)\dx
+\int_\Omega  \big(\vrd \vc{f}\cdot \vud + \vrd G\big) \dx 
 \\ 
&+ \int_\Omega \Big[ \vrd (\vud \otimes \vud):\nabla \vu_0 + p(\vrd,\vtd) \Div  \vu_0 -\tn{S}_\delta(\vtd,\nabla \vud) : \nabla \vu_0 + \vrd\vc{f}\cdot \vu_0  \Big]\dx
\end{aligned} 
\end{equation}
for any $\widetilde{\vt}$ being an extension of the boundary data $\vt_D$ to $\Omega$.

Our aim now is to send $\delta \to 0_+$. We therefore first show estimates independent of $\delta$ and subsequently send $\delta \to 0_+$. We proceed separately for each problem, we present more details for the limit passage in Problem 1, for the other problems we mention only the estimates, as the difficulties with the limit passage are the same as for Problem 1. 

\section{Proof of the main results}

\subsection{Problem 1 (proof of Theorem \ref{t1})}
First, we derive a priori estimates independent of $\delta $, which
will be used afterwards to pass the limit $\delta \rightarrow
0_{+}.$ Note  that due to the form of the pressure and by our construction we have $\varrho_\delta < a$ a.e. in $\Omega$, whence
\begin{equation}
\|\varrho_\delta\|_{\infty} \leq C < \infty.
\end{equation}    
Equation (\ref{weak_energy_1_delta_1}) with the test function $\psi
\equiv 1$ reads%
\begin{equation}
\int_{\Omega }\left( \varrho _{\delta }\mathbf{f}\cdot \mathbf{u}_{\delta
}+\varrho _{\delta }G\right) \,\mathrm{d}x=\int_{\partial \Omega
}L(\vartheta _{\delta }-\vartheta _{0})\,\mathrm{d}S.  \label{energy_1}
\end{equation}%
It implies that%
\begin{equation*}
\int_{\partial \Omega }L\vartheta _{\delta }\,\mathrm{d}S=\int_{\partial
\Omega }L\vartheta _{0}\,\mathrm{d}S+\int_{\Omega }\left( \varrho _{\delta }%
\mathbf{f}\cdot \mathbf{u}_{\delta }+\varrho _{\delta }G\right) \,\mathrm{d}x
\end{equation*}%
which yields%
\begin{equation}
\left\Vert \vartheta _{\delta }\right\Vert _{1,\partial \Omega }\leq C\left(
1+\left\Vert 
\mathbf{u}_{\delta }\right\Vert _{1,p}\right)  \label{esti boundary_1}
\end{equation}%
for some $1\leq p \leq 2$ fixed below.
Using the same test function ($\psi=1$) as in (\ref{weak_entropy_1_delta_1}), we get%
\begin{eqnarray}
&&\int_{\Omega }\left( \frac{\mbox {\F S}_\delta(\vartheta _{\delta },%
\mathbf{u}_{\delta }):\nabla \mathbf{u}_{\delta }}{\vartheta _{\delta }}+ \frac{\Big(
\kappa (\vartheta _{\delta })+\delta (\vartheta _{\delta }^{B}+\vartheta
_{\delta }^{-1})\Big)\left\vert \nabla \vartheta _{\delta }\right\vert ^{2}}{\vartheta_\delta^2}%
+ \frac{\varrho _{\delta }G}{\vartheta _{\delta }}\right) \,\mathrm{d}%
x+\int_{\partial \Omega }\frac{L\vartheta _{0}}{\vartheta _{\delta }}\,%
\mathrm{d}S  \notag \\
&\leq &\int_{\partial \Omega }L\,\mathrm{d}S.  \label{entropy_1 esti}
\end{eqnarray}%
Using Lemma \ref{Korn ineq}, (\ref{visc}), and H\"{o}lder's
inequality, we have%
\begin{equation*}
\left\Vert \mathbf{u}_{\delta }\right\Vert _{1,p}^{p}\leq \left(
\int_{\Omega }\frac{\mbox {\F S}(\vartheta _{\delta },\mathbf{u}%
_{\delta }):\nabla \mathbf{u}_{\delta }}{\vartheta _{\delta }}\,\mathrm{d}%
x\right) ^{\frac{p}{2}}\left( \int_{\Omega }\vartheta _{\delta }^{\frac{%
p(1-\alpha )}{2-p}}\,\mathrm{d}x\right) ^{\frac{2-p}{2}}.
\end{equation*}%
Thus, combining it with the entropy inequality (\ref{entropy_1 esti}) implies
\begin{equation} \label{vel_1}
\left\Vert \mathbf{u}_{\delta }\right\Vert _{1,p}^{p}\leq C \|\vartheta_\delta\|_{\frac{p(1-\alpha )}{2-p}}^{\frac{p(1-\alpha)}{2}}.
\end{equation}
Furthermore, we also deduce from (\ref{entropy_1 esti}) 
\begin{equation*}
\delta \|\vud\|_{1,2}^2 + \left\Vert 
\frac{\varrho _{\delta }G}{\vartheta _{\delta }}\right\Vert _{1}+\left\Vert
\nabla \ln \vartheta _{\delta }\right\Vert _{2}^{2}+\left\Vert \nabla
\vartheta _{\delta }^{\frac{m}{2}}\right\Vert _{2}^{2}+\delta \left(
\left\Vert \nabla \vartheta _{\delta }^{\frac{B}{2}}\right\Vert
_{2}^{2}+\left\Vert \nabla \vartheta _{\delta }^{-\frac{1}{2}}\right\Vert
_{2}^{2}\right) +\left\Vert \vartheta _{\delta }^{-1}\right\Vert
_{1,\partial \Omega }\leq C.
\end{equation*}%
This, together with \eqref{esti boundary_1} implies %
\begin{eqnarray*}
\left\Vert \vartheta _{\delta }\right\Vert _{3m} &\leq &C\left( \left\Vert
\vartheta _{\delta }\right\Vert _{1,\partial \Omega }+\left\Vert \nabla
\vartheta _{\delta }^{\frac{m}{2}}\right\Vert _{2}^{\frac{2}{m}}\right)  \\
&\leq &C\left( 1+\left\Vert \vu_\delta\right\Vert _{1,p}\right)  \\
&\leq &C\left( 1+\left\Vert \vartheta_\delta\right\Vert _{%
\frac{p(1-\alpha)}{2-p}}^{\frac{1-\alpha}{2}}\right) .
\end{eqnarray*}%
Now, take%
\begin{equation*}
3m=\frac{p(1-\alpha )}{2-p}\Rightarrow p=\frac{6m}{3m+1-\alpha }.
\end{equation*}%
Since $p\geq 1$ and $3m\geq 1$, we require that%
\begin{equation*}
m\geq \frac{1}{3}.
\end{equation*}%
It then follows%
\begin{equation*}
\left\Vert \vartheta _{\delta }\right\Vert_{3m} + \left\Vert \vu _{\delta }\right\Vert_{1,p} \leq C.
\end{equation*}%
We control all three unknown quantities; however, we cannot start with the limit passage. The point is that we do not control the pressure. We need to estimate the pressure in some $L^q$-space for $q>1$. To this end, we introduce the Bogovskii
operator. For a given $v \in L^q(\Omega)$, $1<q<\infty$, there exists $\mathcal B(v):=\mathbf{\Phi }$ a solution to the following
Dirichlet boundary value problem 
\begin{eqnarray}
\text{\textrm{div}}\,\mathbf{\Phi }=v-\frac{1}{|\Omega |}
\int_{\Omega }v\,\mathrm{d}x &&\text{ \ in \ }%
\Omega   \notag \\
\mathbf{\Phi }=\mathbf{0} &&\text{ \ on }\partial \Omega 
\label{Bogovskii1}
\end{eqnarray}%
such that%
\begin{equation}
\left\Vert \mathbf{\Phi }\right\Vert _{1,q}\leq C\left\Vert 
v\right\Vert _{q }.
\label{Bogovskii2}
\end{equation}%
We refer, e.g., to \cite[Section 3.3]{Novotnybook} for
the details of the proof. However, we cannot directly use our desired function $\mathcal B(p^\beta(\vrd,\vtd))$ for some $\beta >0$. The point is the lower order term that comes from the second term on the right-hand side of \eqref{Bogovskii1}$_1$. We therefore first use $\mathbf{\Psi } = \mathcal B(\vrd)$, i.e., the solution to
\begin{eqnarray}
\text{\textrm{div}}\,\mathbf{\Psi }=\varrho_\delta-\frac{1}{|\Omega |}%
\int_{\Omega }\varrho_\delta\text{ }\mathrm{d}x &&\text{ \ in \ }%
\Omega   \notag \\
\mathbf{\Psi }=\mathbf{0} &&\text{ \ on }\partial \Omega
\label{Bogovskii3}
\end{eqnarray}%
such that%
\begin{equation}
\left\Vert \mathbf{\Psi }\right\Vert _{1,q}\leq C\left\Vert \varrho_{\delta
}\right\Vert_{q }.
\label{Bogovskii4}
\end{equation}%

We denote $p_\delta := p(\vrd,\vtd)$. As $\varrho_\delta$ is a bounded function, it is not difficult to see that by virtue of the estimates of the velocity and temperature obtained above we have to consider only one difficult term,
$$
\int_\Omega \varrho_\delta p_\delta  \dx \leq C +  \frac{M}{|\Omega|}\int_\Omega p_\delta \dx.
$$
In fact, multiplying both sides of the momentum equation by $\mathbf{\Psi}$ and using (\ref%
{Bogovskii3}) we get
\begin{equation}
\int_{\Omega }\Big[\varrho _{\delta }(\mathbf{u}_{\delta }\otimes \mathbf{u}%
_{\delta }):\nabla \mathbf{\Psi} +p(\varrho _{\delta })\mathrm{div}\,\mathbf{\Psi} -%
\mbox {\F
S}_\delta(\vartheta _{\delta },\nabla \mathbf{u}_{\delta }):\nabla \mathbf{\Psi}
+\varrho _{\delta }\mathbf{f}\cdot \mathbf{\Psi} \Big]\,\mathrm{d}x=0.
\end{equation}%
Then%
\begin{eqnarray*}
&&\int_{\Omega }p _{\delta }\left( \varrho _{\delta }-\frac{1}{%
|\Omega |}\int_{\Omega }\varrho _{\delta }\text{ }\mathrm{d}x\right) dx 
=\int_{\Omega }p _{\delta }\varrho _{\delta }dx-\frac{M}{|\Omega |%
}\int_{\Omega }p _{\delta }dx \\
&=&-\int_{\Omega }\Big[\varrho _{\delta }(\mathbf{u}_{\delta }\otimes 
\mathbf{u}_{\delta }):\nabla \mathbf{\Psi} -\mbox {\F S}_\delta(\vartheta _{\delta
},\nabla \mathbf{u}_{\delta }):\nabla \mathbf{\Psi} +\varrho _{\delta }\mathbf{f}%
\cdot \mathbf{\Psi} \Big]\,\mathrm{d}x. 
\end{eqnarray*}%
Note that%
\begin{eqnarray*}
&&\int_{\Omega }\varrho _{\delta }(\mathbf{u}_{\delta }\otimes 
\mathbf{u}_{\delta }):\nabla \mathbf{\Psi} \mathrm{d}x \\
&\leq &\int_{\Omega }\varrho _{\delta }\left\vert \mathbf{u}_{\delta
}\right\vert ^{2}\left\vert \nabla \mathbf{\Psi} \right\vert \text{ }\mathrm{d}x\leq
\left\Vert \varrho _{\delta }\right\Vert _{\infty }\left\Vert \mathbf{u}%
_{\delta }\right\Vert _{\frac{3p}{3-p}}^{2}\left\Vert \nabla \mathbf{\Psi}
\right\Vert _{\frac{3p}{5p-6}} \\
&\leq &\left\Vert \varrho _{\delta }\right\Vert _{\infty }\left\Vert \mathbf{%
u}_{\delta }\right\Vert _{\frac{3p}{3-p}}^{2}\left\Vert \varrho _{\delta
}\right\Vert _{\frac{3p}{5p-6}}\leq C\left\Vert \mathbf{u}_{\delta
}\right\Vert _{1,p}^{2}\leq C
\end{eqnarray*}%
and%
$$
\begin{aligned}
\int_{\Omega }\mbox {\F S}_\delta (\vartheta _{\delta },\nabla 
\mathbf{u}_{\delta }):\nabla &\mathbf{\Psi} \,\mathrm{d}x\leq C\int_{\Omega
}[(1+\vartheta _{\delta })^{\alpha }+\delta \vartheta_\delta ]|\nabla \mathbf{u}_{\delta }|\text{ }%
\left\vert \nabla \mathbf{\Psi} \right\vert \mathrm{d}x \\
&\leq C\left( 1+\left\Vert \vartheta _{\delta }\right\Vert _{3m}^{\alpha
}\right) \left\Vert \nabla \mathbf{u}_{\delta }\right\Vert _{p}\left\Vert
\varrho _{\delta }\right\Vert _{q} + C \sqrt{\delta} \|\nabla \vud\|_2 \sqrt{\delta}\|\vtd\|_{3B}\Vert\nabla \mathbf{\Psi}\Vert_{q_1}\leq C,
\end{aligned}
$$
where
\begin{equation*}
q=\frac{3mp }{3m(p-1)-\alpha p }.
\end{equation*}%
The denominator is positive provided $m>\frac {1+\alpha}{3}$.
The estimate (\ref{Bogovskii4}) and the properties  of $\mathbf{f}$ and $\varrho _{\delta }$ imply the following bound
\begin{eqnarray*}
\int_{\Omega }\varrho _{\delta }\mathbf{f}\cdot \mathbf{\Psi}\,\mathrm{d}x \leq C.
\end{eqnarray*}

However, due to the condition \eqref{M_limit} we can easily show that the term on the right-hand side above is uniformly bounded. To this aim, we use the argument from \cite{FeNe2018}. Due to condition \eqref{M_limit} we know that there exists $\lambda >1$ such that $\frac{\lambda}{|\Omega|}\int_\Omega \vrd  \dx = \frac{\lambda M}{|\Omega|} <a$. Thus
$$
\begin{aligned}
\int_\Omega \varrho_\delta p_\delta  \dx &\leq C +  \frac{M}{|\Omega|}\int_\Omega p_\delta \dx \\
&= C +  \frac{M}{|\Omega|}\int_{\vrd \leq \lambda \frac{M}{|\Omega|}} p_\delta \dx + \frac{M}{|\Omega|}\int_{\vrd > \lambda \frac{M}{|\Omega|}} p_\delta \dx \\
& \leq C + \frac{\lambda M^2}{|\Omega|^2} h\Big(\frac{\lambda M}{|\Omega|}\Big) \int_\Omega \vtd \dx + \frac 1\lambda \int_\Omega \vrd p_\delta \dx. 
\end{aligned}
$$
This implies that
$$
\int_\Omega \varrho_\delta p_\delta  \dx \leq C,
$$
consequently also
\begin{equation} \label{pressure_1}
\int_\Omega  p_\delta  \dx \leq C.
\end{equation}

We are now ready to use $\mathbf{\Phi } := \mathcal B (p_\delta ^\beta)$ for $\beta >0$  as a test
function in (\ref{weak_mom_delta_1}):%
\begin{eqnarray}
\int_{\Omega }p_{\delta }p_{\delta }^{\beta }\dx
&=&-\int_{\Omega }\varrho _{\delta }(\mathbf{u}_{\delta }\otimes \mathbf{u}%
_{\delta }):\nabla \mathbf{\Phi }\dx+\int_{\Omega }%
\mbox {\F
S}_\delta(\vartheta _{\delta },\nabla\mathbf{u}_{\delta }):\nabla \mathbf{\Phi }%
\dx-\int_{\Omega }\varrho _{\delta }\mathbf{f}\cdot \mathbf{%
\Phi }\dx  \notag \\
&&+\frac{1}{|\Omega |}\left( \int_{\Omega }p_{\delta }^{\beta }\dx\right) \left( \int_{\Omega }p_{\delta }\dx\right)   \notag \\
&=&I_{1}+I_{2}+I_{3}+I_{4}.  \label{moment_1}
\end{eqnarray}%
Then, due to \eqref{pressure_1} the term $I_4$ can be easily controlled by the left-hand side. It is also not difficult  to bound $I_3$ as the sequence of densities is bounded in $L^\infty(\Omega)$. We are therefore left with the first two integrals.
\begin{eqnarray*}
\left\vert I_{1}\right\vert  &\leq &\int_{\Omega }\varrho _{\delta
}\left\vert \mathbf{u}_{\delta }\right\vert ^{2}\left\vert \nabla \mathbf{%
\Phi }\right\vert \text{ }\mathrm{d}x\leq \left\Vert \varrho _{\delta
}\right\Vert _{\infty }\left\Vert \mathbf{u}_{\delta }\right\Vert _{\frac{3p%
}{3-p}}^{2}\left\Vert \nabla \mathbf{\Phi }\right\Vert _{\frac{3p}{5p-6}} \\
&\leq &C\left\Vert \varrho _{\delta }\right\Vert _{\infty }\left\Vert 
\mathbf{u}_{\delta }\right\Vert _{1,p}^{2}\left\Vert p_{\delta }\right\Vert
_{\frac{3p\beta }{5p-6}}^\beta\leq C\left\Vert p_{\delta }\right\Vert _{%
\frac{3p\beta }{5p-6}}^{\beta },
\end{eqnarray*}%
where, in order to get control for some $\beta >0$, we require $\frac{3p\beta}{5p-6}<1+\beta$. This leads to 
\begin{equation} \label{beta_1}
\beta < \frac{5p-6}{6-2p}
\end{equation}
and since $p\leq 2$, the restrictions are $\beta <1$ and
\begin{equation*}
p>\frac{6}{5}\Rightarrow m>\frac{1-\alpha }{2}.
\end{equation*}%
Next
$$
\begin{aligned}
\left\vert I_{2}\right\vert  &\leq C\int_{\Omega }[(1+\vartheta _{\delta
})^{\alpha }+\delta \vartheta_{\delta} ]|\nabla \mathbf{u}_{\delta }|\text{ }\left\vert \nabla \mathbf{\Phi }\right\vert \,\mathrm{d}x\\
&\leq C\left( 1+\left\Vert \vartheta _{\delta
}\right\Vert _{3m}^{\alpha }\right) \left\Vert \nabla \mathbf{u}_{\delta
}\right\Vert _{p}\left\Vert \nabla \mathbf{\Phi }\right\Vert _{q} \\
&+ \delta \|\nabla \vud\|_2 \|\vtd\|_{3B} \|\nabla \mathbf{\Phi}\|_{q_1} \leq C\left\Vert p_{\delta }\right\Vert
_{q\beta }^{\beta },
\end{aligned}
$$
where (recall that $B$ can be taken arbitrarily large)
\begin{equation*}
q=\frac{3mp }{3m(p-1)-\alpha p }.
\end{equation*}%
Thus, we require $\frac{3mp\beta}{3m(p-1)-\alpha p} <1+\beta$ which results in
\begin{equation} \label{beta_2}
\beta < \frac{3m(p-1) -\alpha p}{3m+\alpha p}.
\end{equation}
In order to obtain some $\beta >0$, we therefore require $3m(p-1) -\alpha p>0$, which implies $m > \frac{1+\alpha}{3}$. We are therefore ready for the limit passage in the continuity and momentum equations as well as in the entropy inequality. Recall that we have (for a suitable subsequence which we, however, relabel)
$$
\begin{aligned}
\vu_\delta &\rightharpoonup \vu && \text{ in } W^{1,p}_0(\Omega) \\
\vu_\delta & \to \vu && \text{ in } L^q(\Omega), \, q < \frac{3p}{p-3} \\
\nabla \vartheta_\delta & \rightharpoonup \nabla \vartheta && \text{ in } L^r(\Omega), \, \text{ for some } r>1 \\
\vartheta_\delta &\to \vartheta && \text { in } L^q(\Omega), \, q< 3m \\
\varrho_\delta & \rightharpoonup^* \varrho && \text { in } L^\infty(\Omega) \\
p_\delta & \rightharpoonup \overline{p} && \text{ in } L^{1+\beta}(\Omega)\\
\delta \vtd \nabla \vud &\to 0 && \text{ in } L^1(\Omega) \\ 
\delta \Big(\vtd^B + \frac{1}{\vtd}\Big) \nabla \vtd &\to 0 && \text{ in } L^1(\Omega). 
\end{aligned}
$$
Note that $r=2$ for $m\geq 2$ while $r= \frac{3m}{m+1}$ if $\frac 12 <m<2$.
We find that the limit functions $(\vu,\varrho, \vartheta,\overline{p})$ satisfy
\begin{equation} \label{weak_cont_final1}
\int_\Omega \vr \vu \cdot \nabla \psi \dx = 0
\end{equation}
for all $\psi \in C^1(\overline{\Omega})$;

\begin{equation} \label{weak_mom_final1}
\int_\Omega \Big[ \vr (\vu \otimes \vu):\nabla \vcg{\varphi} + \overline{p} \Div \vcg{\varphi} -\tn{S}(\vt,\nabla \vu) : \nabla \vcg{\varphi} + \vr\vc{f}\cdot \vcg{\varphi}\Big]\dx = 0
\end{equation}
for all $\vcg{\varphi}\in C^1_0(\Omega;\R^3)$;

\begin{equation} \label{weak_entropy_1_final1}
\begin{aligned}
&\int_\Omega \Big(\frac{\tn{S}(\vt,\nabla \vu):\nabla \vu }{\vt} - \frac{\vc{q}(\vt,\nabla \vt)\cdot\nabla \vt }{\vt^2} + \frac{\vr G}{\vt}\Big)\psi \dx \\
&\leq \int_{\partial \Omega} \frac{L(\vt-\vt_0)}{\vt}\psi \dS -\int_\Omega \Big(\frac{\vc{q}(\vt,\nabla \vt)\cdot\nabla \psi}{\vt}+  \overline{\vr s}\vu \cdot \nabla \psi\Big)\dx
\end{aligned} 
\end{equation}
for all non-negative $\psi \in C^1(\overline{\Omega})$, where $\overline{\vr s}$ is the weak limit of $\vrd s(\varrho_\delta, \vartheta_\delta)$. Note that we used the weak lower semicontinuity of several terms in order to pass to the limit on the left-hand side.

\begin{equation} \label{total_energy_integrated_final_1}
\int_{\partial \Omega }L(\vartheta -\vartheta _{0})\,\mathrm{d}%
S=\int_{\Omega }\left( \varrho \mathbf{f}\cdot \mathbf{u}+\varrho G\right) \,%
\mathrm{d}x.
\end{equation}

Indeed, we do not know whether $\overline{p} = p(\varrho,\vartheta)$ and $\overline{\vr s} = \vr s(\varrho,\vartheta)$. Before dealing with this problem, which is equivalent in fact with the strong convergence of the density sequence, let us look at the convergence in the weak formulation of the total energy balance. 

We would like to pass to the limit $\delta \to 0_+$ in the following equality%
\begin{eqnarray}
&&\int_{\Omega }\Big[\varrho _{\delta }\Big(\frac{1}{2}|\mathbf{u}_{\delta
}|^{2}+g(\vartheta _{\delta })\Big)\mathbf{u}_{\delta }
+p(\varrho _{\delta },\vartheta _{\delta })\mathbf{u}_{\delta } -\mbox
{\F S}_\delta(\vartheta _{\delta },\nabla \mathbf{u}_{\delta })\mathbf{u}%
_{\delta } -\kappa _{\delta }(\vartheta _{\delta })\nabla
\vartheta _{\delta }\Big]\cdot \nabla
\psi\dx  \notag \\
&&+\int_{\Omega }\Big(\varrho _{\delta }\mathbf{f}\cdot \mathbf{u}_{\delta
}+\varrho _{\delta }G\Big)\psi \,\mathrm{d}x=\int_{\partial \Omega
}L(\vartheta _{\delta }-\vartheta _{0})\psi \,\mathrm{d}S.
\end{eqnarray}%
Note that the only terms we have to consider in more detail are the second and third terms on the left-hand side. To deal with the latter, we need
$$
\frac{\alpha}{3m} + \frac{1}{p} + \frac{3-p}{3p} <1. 
$$
It easily implies
$$
m>1.
$$
 The term containing the pressure is slightly more complex. Here, we need 
 $$
 \frac{1}{1+\beta} + \frac{3-p}{3p} <1.
 $$
This leads to  $\beta > \frac{3-p}{4p-3}$ and we need to compare this inequality with the inequalities obtained above \eqref{beta_1} and \eqref{beta_2}. We have two inequalities and we need to justify that both of them allow for some $\beta>0$ so that all three inequalities are fulfilled. 

Condition \eqref{beta_1} gives
$$
\frac{5p-6}{6-2p}> \frac{3-p}{4p-3}.
$$
It results in
$$
p>\frac 32, \quad \text{ i.e., } m>1-\alpha.
$$
Condition \eqref{beta_2} gives
$$
\frac{3m(p-1)-\alpha p}{3m+\alpha p}> \frac{3-p}{4p-3}.
$$
This simplifies after plugging the form of $p$ in
$$
\frac{3m-1-\alpha}{3m+1+\alpha} > \frac{m+1-\alpha}{5m-1+\alpha}
$$
which yields $m>1$. Under this assumption, we may pass to the limit in the weak formulation of the total energy balance to get
\begin{eqnarray}
&&\int_{\Omega }\Big[\varrho \Big(\frac{1}{2}|\mathbf{u}|^{2}+g(\vartheta )\Big)\mathbf{u}\cdot \nabla \psi
+\overline{p}\mathbf{u}\cdot \nabla
\psi -\mbox
{\F S}(\vartheta ,\nabla \mathbf{u})\mathbf{u}%
\cdot \nabla \psi -\kappa (\vartheta)\nabla
\vartheta\cdot \nabla \psi \Big]  \notag \\
&&+\int_{\Omega }\Big(\varrho \mathbf{f}\cdot \mathbf{u}
+\varrho G\Big)\psi \,\mathrm{d}x=\int_{\partial \Omega
}L(\vartheta-\vartheta _{0})\psi \,\mathrm{d}S.
\end{eqnarray}%
To conclude, we need to deal with the strong convergence of the density sequence.

However, this task is quite standard nowadays. Exactly as in \cite[Section 7]{CFJP} we deduce the effective viscous flux identity in the form
$$
\overline{ \varrho p } - \varrho \overline{p} = \Big(\frac 43 \mu(\vartheta) + \xi(\vartheta)\Big) (\overline{\varrho \Div \vu} -\varrho \Div \vu). 
$$
It is also easy to see that
$$
\overline{\vartheta \varrho h(\varrho)} = \vartheta \overline{\varrho h(\varrho)}
$$
as well as 
$$
\overline{\vartheta \varrho h(\varrho)\varrho} = \vartheta \overline{\varrho h(\varrho)\varrho}
$$
and exactly as in \cite{CFJP} (see also \cite{Novotnybook}) we can deduce the strong convergence of the density sequence. The proof of Theorem \ref{t1} is complete.

\subsection{Problem 2 (proof of Theorem \ref{t2})}

Now, we turn to Problem 2. We take the ballistic energy inequality (\ref%
{ballistic_energy_1_delta_2}) and choose a particular function $\widetilde{%
\vartheta }:=\vartheta _{L}$ such that (cf. \cite{ChaudAA2022})%
$$
\begin{aligned}
\Delta \vartheta _{L} &=0\text{ \ \ \ \ in }\Omega  \\
\vartheta _{L} &=\vartheta _{D}\text{ \ \ on }\partial \Omega .
\end{aligned}
$$%
As $\vartheta _{D}\in W^{2-\frac{1}{q},q}(\partial \Omega )$ for some $q>3$,
we know that $\vartheta _{L}\in W^{2,q}(\Omega )\hookrightarrow C^{1}(%
\overline{\Omega }).$ Furthermore, by maximum and minimum principles, there
exist $\underline{\theta }$\ $=\min_{x\in \partial \Omega }\vartheta _{D},$
and $\overline{\theta }$\ $=\max_{x\in \partial \Omega }\vartheta _{D}$
such that $\underline{\theta }\leq \vartheta _{L}\leq \overline{\theta }.$
Thus, $\vartheta _{L}$ is an admissible test function in (\ref%
{ballistic_energy_1_delta_2}). We have%
\begin{eqnarray*}
&&\int_{\Omega }\Big( \frac{\mbox {\F S}_\delta(\vartheta _{\delta },%
\nabla \mathbf{u}_{\delta }):\nabla \mathbf{u}_{\delta }}{\vartheta _{\delta }}+%
\frac{\big(\kappa (\vartheta _{\delta })+\delta (\vartheta _{\delta
}^{B}+\vartheta _{\delta }^{-1})\big)\left\vert \nabla \vartheta _{\delta
}\right\vert ^{2}}{\vartheta _{\delta }^{2}}+\frac{\varrho _{\delta }G}{%
\vartheta _{\delta }}\Big) \widetilde{\vartheta }\dx \\
&\geq &\underline{\theta }\int_{\Omega }\Big( \frac{\mbox {\F S}_\delta (\vartheta _{\delta },\mathbf{u}_{\delta }):\nabla \mathbf{u}_{\delta }}{%
\vartheta _{\delta }}+\frac{\big(\kappa (\vartheta _{\delta })+\delta
(\vartheta _{\delta }^{B}+\vartheta _{\delta }^{-1})\big)\left\vert \nabla
\vartheta _{\delta }\right\vert ^{2}}{\vartheta _{\delta }^{2}}+\frac{%
\varrho _{\delta }G}{\vartheta _{\delta }}\Big) \dx \\
&\geq &C\Big(||\mathbf{u}_{\delta }||_{1,p}^{2}||\vartheta _{\delta }||_{\frac{%
p(1-\alpha )}{2-p}}^{-(1-\alpha )}+\big\Vert \nabla \ln \vartheta _{\delta
}\big\Vert _{2}^{2}+\left\Vert \nabla \vartheta _{\delta }^{\frac{m}{2}%
}\right\Vert _{2}^{2}\\
& & +\delta \Big( \big\Vert \nabla \vartheta _{\delta }^{%
\frac{B}{2}}\big\Vert _{2}^{2}+\big\Vert \nabla \vartheta _{\delta }^{-%
\frac{1}{2}}\big\Vert _{2}^{2}+\Vert\mathbf{u}_{\delta}\Vert_{1,2}^2 \Big) +\Big\Vert \frac{\varrho _{\delta }G%
}{\vartheta _{\delta }}\Big\Vert _{1}\Big).
\end{eqnarray*}%
Then%
\begin{eqnarray*}
&&C\left(||\mathbf{u}_{\delta }||_{1,p}^{2}||\vartheta _{\delta }||_{\frac{%
p(1-\alpha )}{2-p}}^{-(1-\alpha )}+\left\Vert \nabla \ln \vartheta _{\delta
}\right\Vert _{2}^{2}+\left\Vert \nabla \vartheta _{\delta }^{\frac{m}{2}%
}\right\Vert _{2}^{2}+\Big\Vert \frac{\varrho _{\delta }G%
}{\vartheta _{\delta }}\Big\Vert _{1}\right)\\ 
&& + C\delta \left( \big\Vert \nabla \vartheta _{\delta }^{%
\frac{B}{2}}\big\Vert _{2}^{2}+\big\Vert \nabla \vartheta _{\delta }^{-%
\frac{1}{2}}\big\Vert _{2}^{2} +\Vert\mathbf{u}_{\delta}\Vert_{1,2}^2 \right) \\
&\leq &\int_{\Omega }\Big( \frac{\big(\kappa (\vartheta _{\delta })+\delta
(\vartheta _{\delta }^{B}+\vartheta _{\delta }^{-1})\big)\nabla \vartheta
_{\delta }\cdot\nabla \vartheta _{L}}{\vartheta _{\delta }}-\varrho _{\delta
}s(\varrho _{\delta },\vartheta _{\delta })\mathbf{u}_{\delta }\cdot \nabla 
\widetilde{\vartheta }+\varrho _{\delta }G+\varrho _{\delta }\mathbf{f\cdot u%
}_{\delta }\Big) \dx.
\end{eqnarray*}%
We estimate the terms on the right-hand side of the above inequality one by
one. We first look at the term with $\kappa (\vartheta _{\delta }).$%
$$
\begin{aligned}
\int_{\Omega }&\frac{\kappa (\vartheta _{\delta })\nabla \vartheta _{\delta
}\cdot\nabla \vartheta _{L}}{\vartheta _{\delta }}\dx =\int_{\Omega }\nabla
K(\vartheta _{\delta })\cdot \nabla \vartheta _{L}\dx \\
&=-\int_{\Omega }K(\vartheta _{\delta }) \Delta \vartheta
_{L}\dx+\int_{\partial \Omega }K(\vartheta _{\delta }) \frac{\partial
\vartheta _{L}}{\partial \mathbf{n}}\,{\rm d}S=\int_{\partial \Omega }K(\vartheta
_{D}) \frac{\partial \vartheta _{L}}{\partial \mathbf{n}}\,{\rm d}S,
\end{aligned}
$$
where%
\begin{equation*}
K^{\prime }(z)=\frac{\kappa (z)}{z}
\end{equation*}%
and the right-hand side is bounded by $C=C(\vartheta _{D},\partial \Omega ).$
Moreover, it is not difficult to find by integrating by parts that%
\begin{equation*}
\int_{\Omega }\Big( \frac{\big(\delta (\vartheta _{\delta }^{B}+\vartheta
_{\delta }^{-1})\big)\nabla \vartheta _{\delta }\cdot\nabla \vartheta _{L}}{%
\vartheta _{\delta }}\Big) \dx= \delta \int_{\partial \Omega} \Big(\frac{\vartheta_D^B}{B}-\frac{1}{\vartheta_D}\Big)\pder{\vartheta_L}{\vc{n}}\,{\rm d}S
\end{equation*}%
and the right-hand side is bounded by $C=C(\vartheta_D,\partial\Omega)$.
Note that the specific entropy $s(\varrho _{\delta },\vartheta _{\delta })$
is given by (\ref{entropy_2}), i.e.,%
\begin{equation*}
s(\varrho ,\vartheta )=\int_{\vartheta _{M}}^{\vartheta }\frac{g^{\prime }(s)%
}{s}\,\mathrm{d}s-\int_{\varrho _{M}}^{\varrho }\frac{h(z)}{z}\,\mathrm{d}z,
\end{equation*}%
where $g(s)$ is given by (\ref{energy_internal_2}). Then, the properties of
the functions $g$ and $h$ imply that%
\begin{equation*}
\varrho s(\varrho ,\vartheta )\leq C\varrho(1+\left[ \ln \varrho \right] ^{+}+\left[ \ln
\vartheta \right] ^{+}),
\end{equation*}%
where $f^{+}=\max \{f(x),0\}.$%
\begin{eqnarray*}
\int_{\Omega }\varrho _{\delta }s(\varrho _{\delta },\vartheta _{\delta })%
\mathbf{u}_{\delta }\cdot \nabla \widetilde{\vartheta }\dx &\leq
&C\int_{\Omega }\varrho _{\delta }(1+\left[ \ln \varrho _{\delta }\right]
^{+}+\left[ \ln \vartheta _{\delta }\right] ^{+})\left\vert \mathbf{u}%
_{\delta }\right\vert |\nabla \widetilde{\vartheta }| \dx \\
&\leq &C\Vert\nabla\widetilde{\vartheta }\Vert_{q}\left\Vert \mathbf{u}_{\delta }\right\Vert _{1,p}+C\int_{\Omega }%
\varrho_\delta\left[ \ln \vartheta _{\delta }\right] ^{+}\left\vert \mathbf{u}_{\delta
}\right\vert\cdot|\nabla \widetilde{\vartheta }| \dx \\
&\leq &C\left\Vert \mathbf{u}_{\delta }\right\Vert _{1,p}+C(\varepsilon)\left\Vert 
\vartheta _{\delta }\right\Vert _{3m}^\varepsilon\left\Vert \mathbf{u}_{\delta
}\right\Vert _{1,p} \\
\end{eqnarray*}%
and%
\begin{equation*}
\int_{\Omega }\varrho _{\delta }\mathbf{f\cdot u}_{\delta }\dx\leq \left\Vert
\varrho _{\delta }\right\Vert _{\infty }\left\Vert \mathbf{u}_{\delta
}\right\Vert _{1,p}.
\end{equation*}%
Therefore, we have%
\begin{eqnarray*}
&&C\Big(||\mathbf{u}_{\delta }||_{1,p}^{2}||\vartheta _{\delta }||_{\frac{%
p(1-\alpha )}{2-p}}^{-(1-\alpha )}+\left\Vert \nabla \ln \vartheta _{\delta
}\right\Vert _{2}^{2}+\big\Vert \nabla \vartheta _{\delta }^{\frac{m}{2}%
}\big\Vert _{2}^{2}\\
&& +\delta \Big( \Big\Vert \nabla \vartheta _{\delta }^{%
\frac{B}{2}}\Big\Vert _{2}^{2}+\Big\Vert \nabla \vartheta _{\delta }^{-%
\frac{1}{2}}\Big\Vert _{2}^{2} + \|\mathbf{u}_{\delta }\|_{1,2}^2 \Big) +\left\Vert \frac{\varrho _{\delta }G%
}{\vartheta _{\delta }}\right\Vert _{1}\Big)\\
&\leq& C\big( 1+\left\Vert \mathbf{u}_{\delta }\right\Vert
_{1,p}+C(\varepsilon )\left\Vert \mathbf{u}_{\delta }\right\Vert _{1,p} \|\vtd\|_{3m}^\varepsilon\big).
\end{eqnarray*}%
We set again $\frac{p(1-\alpha)}{2-p} = 3m$, i.e., $p = \frac{6m}{3m+1-\alpha}$ together with assumption $m> \frac{1-\alpha}{3}$. Then
\begin{equation*}
||\mathbf{u}_{\delta }||_{1,p}^{2}\leq C\big( 1+\left\Vert \mathbf{u}%
_{\delta }\right\Vert _{1,p} ||\vartheta _{\delta }||_{3m}^{1-\alpha +\varepsilon} + \left\Vert \mathbf{u}%
_{\delta }\right\Vert _{1,p}||\vartheta _{\delta }||_{3m}^{1-\alpha}\big).
\end{equation*}%
It follows that%
\begin{equation*}
||\mathbf{u}_{\delta }||_{1,p}\leq C\big( ||\vartheta _{\delta }||_{3m}^{1-\alpha +\varepsilon}+1\big) .
\end{equation*}%
On the other hand, we have%
\begin{eqnarray*}
||\vartheta _{\delta }||_{3m} &\leq &C\Big( \left\Vert \vartheta _{\delta
}\right\Vert _{1,\partial \Omega }+\left\Vert \nabla \vartheta _{\delta }^{%
\frac{m}{2}}\right\Vert _{2}^{\frac{2}{m}}\Big) =C\Big( 1+\left\Vert
\nabla \vartheta _{\delta }^{\frac{m}{2}}\right\Vert _{2}^{\frac{2}{m}%
}\Big)  \\
&\leq &C\big( 1+||\vartheta _{\delta }||_{3m}^{\frac{%
1-\alpha +2\varepsilon}{m}}\big) ,
\end{eqnarray*}%
where we note that on the boundary $\vartheta _{\delta }=\vartheta _{D}$, $\left\Vert
\vartheta _{\delta }\right\Vert _{1,\partial \Omega }$ is bounded by a
constant $C$ that depends on $\vartheta _{D}$ and $\partial \Omega$ and $\varepsilon$ can be taken arbitrarily small, but positive. 
In order to control $||\vartheta _{\delta }||_{3m}$, we require that%
\begin{equation*}
\frac{1-\alpha }{m}<1\Longrightarrow m>1-\alpha .
\end{equation*}%
Then, we have%
\begin{equation*}
||\vartheta _{\delta }||_{3m}\leq C.
\end{equation*}%
Therefore, 
\begin{equation*}
\delta^{\frac 1B} \|\vtd\|_{3B} + \delta\|\vtd^{-3}\|_1^{\frac 13} + \sqrt{\delta} \|\vud\|_{1,2}+ \|\mathbf{u}_{\delta }\|_{1,p}+\|\vartheta _{\delta }\|_{3m}\leq C.
\end{equation*}%
The estimate of pressure is similar to Problem 1, we skip the details. Similarly, the limit passage $\delta \to 0_+$ can be performed as for Problem 1, we first prove the effective viscous flux identity, use compact embedding for the sequence of velocities and temperatures, and based on the fact that the renormalized continuity equations hold for the limit problem, we get the strong convergence of the sequence of densities. The condition allowing the limit passage in the total energy balance are the same as for Problem 1. 

\subsection{Problem 3 (proof of Theorem \ref{t3})}

Next, we deal with Problem 3. We first assume that $0<\alpha<1$. Equation (\ref{weak_energy_3_delta_3}) with the
test function $\psi \equiv 1$ reads%
\begin{eqnarray}
&&\int_{\Omega }\left( \varrho _{\delta }\mathbf{f}\cdot \mathbf{u}_{\delta
}+\varrho _{\delta }G\right) \,\mathrm{d}x-\int_{\partial \Omega
}L(\vartheta _{\delta }-\vartheta _{0})\,\mathrm{d}S \notag \\
&=&\int_{\Omega }\left[
\varrho _{\delta }(\mathbf{u}_{\delta }\otimes \mathbf{u}_{\delta }):\nabla 
\mathbf{u}_{0}-\mbox {\F S}_\delta(\vartheta _{\delta },\nabla \mathbf{u}_{\delta
}):\nabla \mathbf{u}_{0}+\varrho _{\delta }\mathbf{f}\cdot \mathbf{u}_{0}%
\right] \mathrm{d}x.  \label{energy_3}
\end{eqnarray}
This implies that%
\begin{equation}
\begin{aligned} \label{esti problem3}
\int_{\partial \Omega }&L\vartheta _{\delta }\,\mathrm{d}S \leq  C\Big(\|\nabla \vc{u}_0\|_\infty \left\Vert \mathbf{u}_{\delta
}\right\Vert _{1,p}^{2}+\int_{\Omega }\left\vert \mbox {\F S}_\delta%
(\vartheta _{\delta },\nabla \mathbf{u}_{\delta })\right\vert |\nabla \vc{u}_0|\,
\mathrm{d}x+1\Big)  \\
\leq &C\Big(1+\|\nabla \vc{u}_0\|_\infty \left\Vert \mathbf{u}_{\delta }\right\Vert _{1,p}^{2}+\|\nabla \vc{u}_0\|_\infty \left\Vert 
\mathbf{u}_{\delta }\right\Vert _{1,p}\left\Vert \vartheta _{\delta
}\right\Vert _{\frac{p(1-\alpha )}{2-p}}^{\alpha } + 
\delta \|\nabla \vu_0\|_\infty \|\nabla \vud\|_2 \|\vartheta _{\delta }\|_2 \Big),  
\end{aligned}
\end{equation}
where we note that for $\alpha \in  (0,1)$, $1<p<2$ and $m> \frac {1+\alpha}{3}$ (i.e., $p>1+\alpha$)%
\begin{equation*}
\frac{1}{p}+\frac{\alpha (2-p)}{p(1-\alpha )}=\frac{1+\alpha -p\alpha }{%
p(1-\alpha )}<\frac{p-p\alpha }{p(1-\alpha )}=1.
\end{equation*}%
Using the same test function $\psi \equiv 1$ in (\ref{weak_entropy_1_delta_3}) gives
\begin{eqnarray} \label{estigeneralalpha}
&&\int_{\Omega }\Big( \frac{\mbox {\F S}_{\delta }(\vartheta _{\delta },%
\mathbf{u}_{\delta }):\nabla \mathbf{u}_{\delta }}{\vartheta _{\delta }}+%
\frac{\big(\kappa (\vartheta _{\delta })+\delta (\vartheta _{\delta
}^{B}+\vartheta _{\delta }^{-1})\big)\left\vert \nabla \vartheta _{\delta
}\right\vert ^{2}}{\vartheta _{\delta }^{2}}+\frac{\varrho _{\delta }G}{%
\vartheta _{\delta }}\Big) \,\mathrm{d}x+\int_{\partial \Omega }\frac{%
L\vartheta _{0}}{\vartheta _{\delta }}\,\mathrm{d}S  \notag \\
&\leq &\int_{\partial \Omega }L\,\mathrm{d}S.
\end{eqnarray}%
Then%
\begin{equation*}
\delta \|\vud\|_{1,2}^2 +\left\Vert \frac{\varrho _{\delta }G}{\vartheta _{\delta }}\right\Vert
_{1}+\left\Vert \nabla \ln \vartheta _{\delta }\right\Vert
_{2}^{2}+\Big\Vert \nabla \vartheta _{\delta }^{\frac{m}{2}}\Big\Vert
_{2}^{2}+\delta \Big( \Big\Vert \nabla \vartheta _{\delta }^{\frac{B}{2}%
}\Big\Vert _{2}^{2}+\Big\Vert \nabla \vartheta _{\delta }^{-\frac{1}{2}%
}\Big\Vert _{2}^{2}\Big) +\left\Vert \vartheta _{\delta
}^{-1}\right\Vert _{1,\partial \Omega }\leq C.
\end{equation*}
The only problematic term on the right-hand side of \eqref{esti problem3} is the last one. Using \eqref{estigeneralalpha} we, however, have
$$
\begin{aligned}
C\delta \|\nabla u_0\|_\infty \|\nabla \vud\|_2 \|\vtd\|_2 &\leq C \sqrt{\delta} \|\vud\|_2 \sqrt{\delta} \|\vtd\|_{3B} \leq C \sqrt{\delta} \Big(\|\vtd\|_{1,\partial\Omega} + \|\nabla \vtd^{\frac B2}\|_2^{\frac 2B}\Big) \\
&\leq C \sqrt{\delta} \|\vtd\|_{1,\partial\Omega} + C.
\end{aligned}
$$
We plug this estimate in \eqref{esti problem3} and for sufficiently small $\delta$ we transfer the first term on the right-hand side above to the left-hand side of \eqref{esti problem3}.
Similarly as in Problem 1, we have due to \eqref{estigeneralalpha}%
\begin{equation}\label{79}
\left\Vert \mathbf{u}_{\delta }\right\Vert _{1,p}^{p}\leq C\Vert \vartheta
_{\delta }\Vert _{\frac{p(1-\alpha )}{2-p}}^{\frac{p(1-\alpha )}{2}}.
\end{equation}%
Setting again $\frac{p(1-\alpha)}{2-p}=3m$, i.e., $p=\frac{6m}{3m+1-\alpha}$, and assuming $m> \frac{1-\alpha}{3}$ (but recall that we assume $m >\frac{1+\alpha}{3}$) we get, again combining \eqref{79} with \eqref{esti problem3}
\begin{eqnarray*}
\left\Vert \vartheta _{\delta }\right\Vert _{3m} &\leq &C\Big( \left\Vert
\vartheta _{\delta }\right\Vert _{1,\partial \Omega }+\Big\Vert \nabla
\vartheta _{\delta }^{\frac{m}{2}}\Big\Vert _{2}^{\frac{2}{m}}\Big)  \\
&\leq &C\Big( 1+\left\Vert \mathbf{u}_{\delta }\right\Vert
_{1,p}^{2}+\left\Vert \mathbf{u}_{\delta }\right\Vert _{1,p}\left\Vert
\vartheta _{\delta }\right\Vert _{3m}^{\alpha }\Big) 
\\
&\leq &C\Big( 1+\left\Vert \vartheta _{\delta }\right\Vert _{3m}^{1-\alpha }+\left\Vert \vartheta _{\delta }\right\Vert _{%
3m}^{\frac{1+\alpha }{2}}\Big) .
\end{eqnarray*}%
Then, it follows that%
\begin{equation*}
\left\Vert \vartheta _{\delta }\right\Vert _{3m}+\left\Vert \mathbf{u}%
_{\delta }\right\Vert _{1,p}\leq C.
\end{equation*}%
The remaining part is identical to Problem 1.

Clearly, the case $\alpha =0$ cannot be estimated in the same way. The problematic term is $\|\vtd\|_{3m}^{1-\alpha}|_{\alpha=0}$ above. We therefore follow the idea from \cite{CFJP}. We recall Lemma A1 from \cite{CF} which we reformulate to our situation.
\begin{lem}\label{A1}
Let $\Omega \subset \R^3$ be a bounded Lipschitz domain. Let $\vc{U} \in W^{1,p}$, $p>3$ be given so that $\vc{U}\cdot\vc{n}=0$ on $\partial \Omega$. Let $q\in (1,\infty)$. Then for any $\omega >0$ there exists $\vc{u}_\omega \in W^{1,p}(\Omega)$ such that
\begin{itemize}
\item $\vc{u}_\omega = \vc{U}$ on $\partial \Omega$ and $\Div \vc{u}_\omega = 0$ in $\Omega$
\item $\|\vc{u}_\omega\|_q \leq \omega$ and $\|\vc{u}_\omega\|_{1,p} \leq C \|\vc{U}\|_{1,p}$.
\end{itemize}
\end{lem}
We now take $\omega$ sufficiently small and replace our $\vc{u}_0$ with $\vc{u}_\omega$. We return to \eqref{energy_3} and replace the convective term with
$$
\int_\Omega -\vrd (\vud\cdot \nabla \vud) \cdot \vc{u}_\omega \dx,
$$
where we used the weak formulation of the continuity equation. We need at this stage that the term above is integrable, which results in (recall that $\alpha =0$  implies $p=\frac{6m}{3m+1}$)
$$
\frac{3m+1}{6m} + \frac{m+1}{6m} <1,
$$
i.e., $m>1$. Then
$$
\Big|\int_\Omega \vrd (\vud\cdot \nabla \vud) \cdot \vc{u}_\omega \dx\Big| \leq \|\vrd\|_\infty 
\|\vud\|_{1,p}^2 \|\vu_\omega\|_q
$$
for some $q$ sufficiently large. We use this estimate in \eqref{esti problem3} and then repeat the estimates for the temperature
\begin{eqnarray*}%
\left\Vert \vartheta _{\delta }\right\Vert _{3m} &\leq &C\Big( \left\Vert
\vartheta _{\delta }\right\Vert _{1,\partial \Omega }+\Big\Vert \nabla
\vartheta _{\delta }^{\frac{m}{2}}\Big\Vert _{2}^{\frac{2}{m}}\Big)  \\
&\leq &C\Big( 1+\omega \left\Vert \mathbf{u}_{\delta }\right\Vert
_{1,p}^{2}+\left\Vert \mathbf{u}_{\delta }\right\Vert _{1,p}\Big) 
\\
&\leq &C\Big( 1+\omega\left\Vert \vartheta _{\delta }\right\Vert _{3m}+\left\Vert \vartheta _{\delta }\right\Vert _{%
3m}^{\frac{1}{2}}\Big) .
\end{eqnarray*}%
Finally, we fix $\omega$ sufficiently small so that $C\omega <\frac 12$. Thus, the estimate of the temperature follows and consequently, also the estimate of the velocity.
The rest of the proof is identical to the case $0<\alpha<1$, but we get a stronger restriction $m>1$, which also ensures the limit passage in the total energy balance.

For $\alpha =1$ ($p=2$) we combine \eqref{esti problem3} with the estimate $\|\vud\|_{1,2} \leq C$ that follows from \eqref{estigeneralalpha}. We get
$$
\begin{aligned}
\|\vtd\|_{3m} &\leq C \Big(\|\vtd\|_{1,\partial \Omega} + \big\|\nabla\vtd^{\frac m2}\big\|_2^{\frac 2m}\Big) \\
& \leq C\big( 1 + \|\nabla \vu_0\|_\infty \|\vud\|_{1,2}^2 + \|\nabla \vu_0\|_\infty \|\vud\|_{1,2} \|\vtd\|_2\big)
\end{aligned}
$$
and for $m \geq \frac 23$ and $\|\nabla \vu_0\|_\infty$ sufficiently small, we get
$$
\|\vtd\|_{3m} \leq C.
$$
The rest of the proof is similar to the proof of Theorem 1. For $m>1$ we can perform the limit passage in the total energy balance.
Thus, we complete the proof of
Theorem \ref{t3}.

\subsection{Problem 4 (proof of Theorem \ref{t4})}

Finally, we deal with Problem 4. We take the ballistic energy inequality (%
\ref{ballistic_energy_2_delta_4}) and choose a particular function $%
\widetilde{\vartheta }:=\vartheta _{L}$ exactly as in Problem 2. Whence $\vartheta _{L}$ is an admissible test function in (\ref%
{ballistic_energy_2_delta_4}). We also choose $\alpha =1$, i.e., we take $p=2$. We have (recall that the regularization of the form $\tn{S}_\delta$ is meaningless in this case)%
\begin{eqnarray*}
&&\int_{\Omega }\Big( \frac{\mbox {\F S}(\vartheta _{\delta },%
\nabla \mathbf{u}_{\delta }):\nabla \mathbf{u}_{\delta }}{\vartheta _{\delta }}+%
\frac{\big(\kappa (\vartheta _{\delta })+\delta (\vartheta _{\delta
}^{B}+\vartheta _{\delta }^{-1})\big)\left\vert \nabla \vartheta _{\delta
}\right\vert ^{2}}{\vartheta _{\delta }^{2}}+\frac{\varrho _{\delta }G}{%
\vartheta _{\delta }}\Big) \widetilde{\vartheta }\dx \\
&\geq &\underline{\theta }\int_{\Omega }\Big( \frac{\mbox {\F S}(\vartheta _{\delta },\nabla \mathbf{u}_{\delta }):\nabla \mathbf{u}_{\delta }}{%
\vartheta _{\delta }}+\frac{\big(\kappa (\vartheta _{\delta })+\delta
(\vartheta _{\delta }^{B}+\vartheta _{\delta }^{-1})\big)\left\vert \nabla
\vartheta _{\delta }\right\vert ^{2}}{\vartheta _{\delta }^{2}}+\frac{%
\varrho _{\delta }G}{\vartheta _{\delta }}\Big) \dx \\
&\geq &C\Big(||\mathbf{u}_{\delta }||_{1,2}^{2}+\left\Vert \nabla \ln \vartheta
_{\delta }\right\Vert _{2}^{2}+\big\Vert \nabla \vartheta _{\delta }^{%
\frac{m}{2}}\big\Vert _{2}^{2}+\delta \Big( \big\Vert \nabla \vartheta
_{\delta }^{\frac{B}{2}}\big\Vert _{2}^{2}+\big\Vert \nabla \vartheta
_{\delta }^{-\frac{1}{2}}\big\Vert _{2}^{2}\Big) +\Big\Vert \frac{%
\varrho _{\delta }G}{\vartheta _{\delta }}\Big\Vert _{1}\Big).
\end{eqnarray*}%
Furthermore
\begin{eqnarray*}
&&\int_{\Omega }\left[ \varrho _{\delta }(\mathbf{u}_{\delta }\otimes 
\mathbf{u}_{\delta }):\nabla \mathbf{u}_{0}+p(\varrho _{\delta },\vartheta
_{\delta })\mathrm{div}\,\mathbf{u}_{0}-\mbox {\F S}(\vartheta _{\delta
},\nabla \mathbf{u}_{\delta }):\nabla \mathbf{u}_{0}+\varrho _{\delta }%
\mathbf{f}\cdot \mathbf{u}_{0}\right] \mathrm{d}x \\
&\leq & \|\vrd\|_\infty \left\Vert \nabla \mathbf{u}_{0}\right\Vert _{\infty }\left\Vert 
\mathbf{u}_{\delta }\right\Vert _{2}^{2}+\left\Vert \mathrm{div}\,\mathbf{u}%
_{0}\right\Vert _{\infty }\left\Vert p_\delta\right\Vert _{1}+\left\Vert \nabla \mathbf{u}_{0}\right\Vert _{\infty
}\left\Vert \mbox {\F S}(\vartheta _{\delta },\nabla \mathbf{u}_{\delta
})\right\Vert _{1}+C.
\end{eqnarray*}%
Then%
\begin{eqnarray*}
&&C\Big(||\mathbf{u}_{\delta }||_{1,2}^{2}+\big\Vert \nabla \ln \vartheta
_{\delta }\big\Vert _{2}^{2}+C\left\Vert \nabla \vartheta _{\delta }^{%
\frac{m}{2}}\right\Vert _{2}^{2}+\delta \big( \big\Vert \nabla \vartheta
_{\delta }^{\frac{B}{2}}\big\Vert _{2}^{2}+\big\Vert \nabla \vartheta
_{\delta }^{-\frac{1}{2}}\big\Vert _{2}^{2}\big) +\left\Vert \frac{%
\varrho _{\delta }G}{\vartheta _{\delta }}\right\Vert _{1}\Big) \\
&\leq &\int_{\Omega }\Big( \frac{\big(\kappa (\vartheta _{\delta })+\delta
(\vartheta _{\delta }^{B}+\vartheta _{\delta }^{-1})\big)\nabla \vartheta
_{\delta }\cdot\nabla \vartheta _{L}}{\vartheta _{\delta }}+\varrho _{\delta
}s(\varrho _{\delta },\vartheta _{\delta })\mathbf{u}_{\delta }\cdot \nabla 
\vartheta_L+\varrho _{\delta }G+\varrho _{\delta }\mathbf{f\cdot u%
}_{\delta }\Big) \dx \\
&&+C_1 \Big(\|\vrd\|_\infty
\left\Vert \nabla \mathbf{u}_{0}\right\Vert _{\infty }\left\Vert \mathbf{%
u}_{\delta }\right\Vert _{1,2}^{2}+\left\Vert \mathrm{div}\,\mathbf{u}%
_{0}\right\Vert _{\infty }\left\Vert p_\delta\right\Vert _{1}+\left\Vert \nabla \mathbf{u}_{0}\right\Vert _{\infty
}\left\Vert \nabla \mathbf{u}_{\delta }\right\Vert _{2}\left\Vert \vartheta
_{\delta }\right\Vert _{2}+1\Big).
\end{eqnarray*}%
Next, we look at the term with $\kappa (\vartheta _{\delta })$ in the
integrals.
\begin{eqnarray*}
\int_{\Omega }\frac{\kappa (\vartheta _{\delta })\nabla \vartheta _{\delta
}\cdot\nabla \vartheta _{L}}{\vartheta _{\delta }}\dx &=&\int_{\Omega }\nabla
K(\vartheta _{\delta }) \cdot \nabla \vartheta _{L}\dx \\
&=&-\int_{\Omega }K(\vartheta _{\delta }) \Delta \vartheta
_{L}\dx+\int_{\partial \Omega }K(\vartheta _{\delta }) \frac{\partial
\vartheta _{L}}{\partial \mathbf{n}}\,{\rm d}S=\int_{\partial \Omega }K(\vartheta
_{D}) \frac{\partial \vartheta _{L}}{\partial \mathbf{n}}\,{\rm d}S,
\end{eqnarray*}%
where%
\begin{equation*}
K^{\prime }(z)=\frac{\kappa (z)}{z},
\end{equation*}%
and the right-hand side is bounded by $C=C(\vartheta _{D},\partial \Omega ).$
Moreover, it is not difficult to find by integrating by parts that%
\begin{equation*}
\int_{\Omega }\Big( \frac{\big(\delta (\vartheta _{\delta }^{B}+\vartheta
_{\delta }^{-1})\big)\nabla \vartheta _{\delta }\cdot\nabla \vartheta _{L}}{%
\vartheta _{\delta }}\Big) \dx=\delta \int_{\partial \Omega} \Big(\frac{\vartheta_D^B}{B}-\frac{1}{\vartheta_D}\Big)\pder{\vartheta_L}{\vc{n}}\,{\rm d}S.
\end{equation*}%
Note that the specific entropy $s(\varrho _{\delta },\vartheta _{\delta })$
is given by (\ref{entropy_2}), i.e.,%
\begin{equation*}
s(\varrho ,\vartheta )=\int_{\vartheta _{M}}^{\vartheta }\frac{g^{\prime }(s)%
}{s}\,\mathrm{d}s-\int_{\varrho _{M}}^{\varrho }\frac{h(z)}{z}\,\mathrm{d}z,
\end{equation*}%
where $g(s)$ is given by (\ref{energy_internal_2}). The properties of
the functions $g$ and $h$ imply that%
\begin{equation*}
\varrho s(\varrho ,\vartheta )\leq C\varrho (1+\left[ \ln \varrho \right] ^{+}+\left[ \ln
\vartheta \right] ^{+}).
\end{equation*}%
This yields
\begin{eqnarray*}
\int_{\Omega }\varrho _{\delta }s(\varrho _{\delta },\vartheta _{\delta })%
\mathbf{u}_{\delta }\cdot \nabla \vartheta_L\dx &\leq
&C\int_{\Omega }\varrho _{\delta }(1+\left[ \ln \varrho _{\delta }\right]
^{+}+\left[ \ln \vartheta _{\delta }\right] ^{+})\left\vert \mathbf{u}%
_{\delta }\right\vert \dx \\
&\leq &C\left\Vert \mathbf{u}_{\delta }\right\Vert _{1,2}+C\int_{\Omega }%
\left[ \ln \vartheta _{\delta }\right] ^{+}\left\vert \mathbf{u}_{\delta
}\right\vert\dx \\
&\leq &C\left\Vert \mathbf{u}_{\delta }\right\Vert _{1,2}+C(\varepsilon)\left\Vert
\vartheta _{\delta }\right\Vert _{3m}^\varepsilon\left\Vert \mathbf{u}_{\delta
}\right\Vert _{1,2} 
\end{eqnarray*}%
and%
\begin{equation*}
\int_{\Omega }\varrho _{\delta }\mathbf{f\cdot u}_{\delta }\dx\leq \left\Vert
\varrho _{\delta }\right\Vert _{\infty }\left\Vert \mathbf{u}_{\delta
}\right\Vert _{1,2}.
\end{equation*}%
Therefore, we have%
\begin{equation} \label{estimate4}
\begin{aligned}
&C\Big(||\mathbf{u}_{\delta }||_{1,2}^{2}+\Vert \nabla \ln \vartheta
_{\delta }\Vert _{2}^{2}+C\big\Vert \nabla \vartheta _{\delta }^{%
\frac{m}{2}}\big\Vert _{2}^{2}+\delta \Big( \big\Vert \nabla \vartheta
_{\delta }^{\frac{B}{2}}\big\Vert _{2}^{2}+\big\Vert \nabla \vartheta
_{\delta }^{-\frac{1}{2}}\big\Vert _{2}^{2}\Big) +\left\Vert \frac{%
\varrho _{\delta }G}{\vartheta _{\delta }}\right\Vert _{1}\Big)  \\
&\leq C\Big(\left\Vert \mathbf{u}_{\delta }\right\Vert _{1,2}+\|\vrd\|_\infty\left\Vert \nabla 
\mathbf{u}_{0}\right\Vert _{\infty }\left\Vert \mathbf{u}_{\delta
}\right\Vert _{1,2}^{2}+\left\Vert \mathrm{div}\,\mathbf{u}_{0}\right\Vert
_{\infty }\left\Vert p_\delta\right\Vert
_{1}\\
&+\left\Vert \nabla \mathbf{u}_{0}\right\Vert _{\infty }\left\Vert \nabla 
\mathbf{u}_{\delta }\right\Vert _{2}\left\Vert \vartheta _{\delta
}\right\Vert _{2}+\|\vtd\|_{3m}^{\varepsilon} \|\vud\|_{1,2} + \|\vud\|_{1,2} + 1 \Big).  
\end{aligned}
\end{equation}
The control of pressure is similar to Problem 1. Using $\mathbf{\Phi }= \mathcal B(\vrd)$  as a test function in (\ref{weak_mom_delta_4}) gives
\begin{eqnarray*}
\int_{\Omega }p_{\delta }\vrd\dx
&=&-\int_{\Omega }\varrho _{\delta }(\mathbf{u}_{\delta }\otimes \mathbf{u}%
_{\delta }):\nabla \mathbf{\Phi }\dx+\int_{\Omega }%
\mbox {\F
		S}(\vartheta _{\delta },\nabla \mathbf{u}_{\delta }):\nabla \mathbf{%
\Phi }\dx-\int_{\Omega }\varrho _{\delta }\mathbf{f}\cdot 
\mathbf{\Phi }\dx  \notag \\
&&+\frac{M}{|\Omega |}\left( \int_{\Omega }p_{\delta }^{\beta }\text{ }%
\mathrm{d}x\right) \dx.   \\
\end{eqnarray*}%
Exactly as in the proof of Theorem \ref{t1}, we may show that
\begin{equation} \label{84c}
\int_{\Omega }p_{\delta }\vrd\dx + \int_{\Omega }p_{\delta }\dx \leq C\big(\|\vrd\|_\infty \|\vud\|_{1,2}^2 + \|\vud\|_{1,2}\|\vtd\|_2 + \|\vrd\|_\infty\big).
\end{equation}
We employ the estimate \eqref{84c} and get that if $\|\nabla \vu_0\|_\infty$ is sufficiently small and $m\geq 2$ (this comes from the estimate of the term $\|\vud\|_{1,2} \|\vtd\|_2$) we can estimate the left-hand side of \eqref{estimate4} and thus also the left-hand side of \eqref{84c}.

Finally, in (\ref{weak_mom_delta_4}) we use the test function $\mathcal B(p_\delta^\beta)$ and get

\begin{eqnarray}
\int_{\Omega }p_{\delta }p_{\delta }^{\beta }\text{ }\mathrm{d}x
&=&-\int_{\Omega }\varrho _{\delta }(\mathbf{u}_{\delta }\otimes \mathbf{u}%
_{\delta }):\nabla \mathbf{\Phi }\text{ }\mathrm{d}x+\int_{\Omega }%
\mbox {\F
		S}(\vartheta _{\delta },\mathbf{u}_{\delta }):\nabla \mathbf{%
\Phi }\text{ }\mathrm{d}x-\int_{\Omega }\varrho _{\delta }\mathbf{f}\cdot 
\mathbf{\ \Phi }\,\mathrm{d}x  \notag \\
&&+\frac{1}{|\Omega |}\left( \int_{\Omega }p_{\delta }^{\beta }\text{ }%
\mathrm{d}x\right) \left( \int_{\Omega }p_{\delta }\text{ }\mathrm{d}%
x\right)   \notag \\
&=&I_{1}+I_{2}+I_{3}+I_{4}.
\end{eqnarray}%

We estimate these terms to control the pressure in a better space than $L^1$, then we pass to the limit as in previous results.

It is sufficient to consider the first two integrals. Therefore, 
\begin{eqnarray*}
\left\vert I_{1}\right\vert  &\leq &\int_{\Omega }\varrho _{\delta
}\left\vert \mathbf{u}_{\delta }\right\vert ^{2}\left\vert \nabla \mathbf{%
\Phi }\right\vert \, \mathrm{d}x\leq \left\Vert \varrho _{\delta
}\right\Vert _{\infty }\left\Vert \mathbf{u}_{\delta }\right\Vert _{6}^{2}\left\Vert \nabla \mathbf{\Phi }\right\Vert _{\frac 32} \\
&\leq &C\left\Vert \varrho _{\delta }\right\Vert _{\infty }\left\Vert 
\mathbf{u}_{\delta }\right\Vert _{1,p}^{2}\left\Vert p_{\delta }\right\Vert
_{\frac 32 \beta}^{\beta },
\end{eqnarray*}%
where, in order to get control for some $\beta >0$, we require $\frac 32\beta <1+\beta $. This leads to 
\begin{equation}
\beta <2.
\end{equation}%
Next 
\begin{eqnarray*}
\left\vert I_{2}\right\vert  &\leq &C\int_{\Omega }(1+\vartheta _{\delta
})|\nabla \mathbf{u}_{\delta }|\, \left\vert \nabla \mathbf{\Phi }%
\right\vert \, \mathrm{d}x\leq C\left( 1+\left\Vert \vartheta _{\delta
}\right\Vert _{3m}\right) \left\Vert \nabla \mathbf{u}_{\delta }\right\Vert
_{2}\left\Vert \nabla \mathbf{\Phi }\right\Vert _{\frac{6m}{3m-2}} \\
&\leq &C\left( 1+\left\Vert \vartheta _{\delta }\right\Vert _{3m}\right)
\left\Vert \nabla \mathbf{u}_{\delta }\right\Vert _{2}\left\Vert p_{\delta
}\right\Vert _{\frac{6m\beta}{3m-2} }^{\beta }.
\end{eqnarray*}%
Thus, we require $\frac{3m\beta }{3m-2}<1+\beta $ which results in 
\begin{equation}
\beta <\frac{3m-2}{3m+2}.
\end{equation}%
In order to get some $\beta >0$ we therefore need that $m>\frac{2}{3}$; note, however, that we already have $m\geq 2$. Then, we obtain%
\begin{equation*}
\left\Vert p_{\delta }\right\Vert _{1+\beta }\leq C\Big( \left\Vert \mathbf{%
u}_{\delta }\right\Vert _{1,p}^{2}+\left\Vert \mathbf{u}_{\delta
}\right\Vert _{1,p}\left\Vert \vartheta _{\delta }\right\Vert _{3m} +1\Big) \leq C.
\end{equation*}%
Now, we can start with the limit passage procedure similarly as in Problem
1. Note that for $m\geq 2$ the limit passage in the total energy balance does not cause any troubles.
\section{Acknowledgments}
Zhengguang Guo was funded by Basic Research Program of Jiangsu under Grant No. \linebreak BK20251938. The work of Milan Pokorn\'y was partially supported by the Czech Science Foundation, project No. 25-16592S.

\end{document}